\documentstyle[]{article}
\newcommand{\lemma}[2]{{\bf Lemma #1.}\hskip0.1cm{\it#2}}
\newcommand{\theorem}[2]{{\bf Theorem #1.}\hskip0.1cm{\it#2}}

\newcommand{\remark}[1]{{\bf Remark #1.}\hskip0.3cm}
\newcommand{\proof}{{\bf Proof.}\hskip0.2cm}
\newcommand{\pf}[1]{{\bf Proof of {#1}}.\hskip0.3cm}
\newcommand{\corollary}[2]{{\bf Corollary #1.}\hskip
        0.1cm{\it#2}}
\newcommand{\fd}[2]{\frac{\displaystyle #1}{\displaystyle #2}}

\newcommand{\mb}[1]{\hsa\mbox{#1}\hsa}
\newcommand{\hsa}{\hskip0.19cm}

\newcommand{\cc}{{\cal C}}

\newcommand{\dd}[1]{|\nabla|^{#1}}

\title{
\Large \bf Existence of Smooth Solutions of the Navier-Stokes
Equations}

\author {
 \large Dongsheng Li\thanks{
 Supported by NSF of China: 11171266} \\
 {\small\it Department of mathematics, Xi'an Jiaotong University}\\
 {\small\it Xi'an 710049, China}\\
 {\small\it E-mail:\quad lidsh@mail.xjtu.edu.cn}
 }

\date{}

\begin{document}

\maketitle

\begin{center}
\parbox{12cm}{\small
{{\bf Abstract.} In this paper, we prove existence of smooth
solutions of the Navier-Stokes equations that gives a positive answer to the problem
proposed by Fefferman [3].

\bigskip


}}
\end{center}

\section{Introduction and the main results}

The Navier-Stokes equations are given by
$$\left\{\begin{array}{l}u_t-\nu\Delta u+(u\cdot \nabla)u+\nabla p=f;\\
\mbox{div}u=0,\end{array}\right.\eqno{(1.1)}$$ where $\nu$ is a
positive constant. The existence of smooth solutions of (1.1) is an
open problem standing for a long time. Here we only mention some
remarkable works that tried to solve it. Leray [4] showed existence
of weak solutions. Under additional assumptions of more
integrability of $u$, Serrin [7] proved existence of smooth
solutions. In [2], Constantin-Fefferman showed smoothness of
solutions with a constraint on vorticity. Caffarelli-Kohn-Nirenberg
[1] gave a partial regularity result of the Navier-Stokes equations
that the dimension of the set of singular points is at most one,
which improved the results of Scheffer [6]. Later Lin [5] simplified
the proof. In this paper, we will solve this long standing problem. We state our main results as the following, which are
corresponding to the statements (A) and (B) in [3] respectively.

\bigskip

\theorem{1.1}{Let $u_0$ be any smooth, divergence-free vector field
in $R^3$ satisfying
$$|\partial^\alpha_xu_0(x)|\leq C_{\alpha K}(1+|x|)^{-K}\mb{on}R^3,\mb{for any}\alpha\mb{and}K.\eqno{(1.2)}$$ Then
there exist smooth functions $p(x,t)$ and $u(x,t)$ on
$R^3\times[0,\infty)$ that satisfy (1.1) with $u(x,0)=u_0(x)$,
$f\equiv0$ and
$$\int_{R^3}|u(x,t)|^2dx<C\mb{for all}t\geq0.\eqno{(1.3)}$$}

\bigskip

\theorem{1.2}{Let $u_0$ be any smooth, divergence-free vector field
in $R^3$ satisfying $$u_0(x+e_j)=u_0(x)\mb{for}1\leq j\leq
3.\eqno{(1.4)}$$ Then there exist smooth functions $p(x,t)$ and
$u(x,t)$ on $R^3\times[0,\infty)$ that satisfy (1.1) with
$u(x,0)=u_0(x)$, $f\equiv0$ and
$$u(x,t)=u(x+e_j,t)\mb{on}R^3\times[0,\infty)\mb{for}1\leq j\leq 3.\eqno{(1.5)}$$}

\bigskip

\remark{1.3} (i) To prove Theorem 1.1 and 1.2, it is the key to get an a priori estimate with sufficient regularity.
In order to express our idea clearly and neatly, we
only focus on solving the problems proposed in [3] which are essential
as considering this kind of a priori estimates of
regularity of the Navier-Stokes equations. However,
using our method, it is not hard to obtain the regularity of
Navier-Stokes equations in higher spatial dimensions including
interior estimates and boundary estimates, and the regularity of
steady-states, where the righthand term $f$ need not to be vanishing.

(ii) In this paper, we only prove Theorem 1.1, while Theorem 1.2 can
be proved similarly. $\Box$

\bigskip

The difficulty of proving the existence of smooth solutions of the
Navier-Stokes equations arises from the following fact: When we
multiply $F$ on both sides of the equation, where $F$ may
contain $u$ and (or) derivative of $u$, and then integrate over
$R^3$, the {\it bad terms} coming from $(u\cdot\nabla)u$ and $\nabla p$ can not be
controlled by the {\it good terms} coming from $u_t$ and $-\Delta u$
(except $F=u$).

To overcome this difficulty, we multiply a series of $\{F_k\}$ to the equation and take integral over $R^3$.
Then we have infinitely many inequalities and the {\it bad terms} in the former
inequalities can be controlled by the {\it good terms} in the
later inequalities. Therefore if we add all the inequalities together, all the {\it bad terms} can be
controlled. To dealt with the nonlinear term with differential operators, it is convenient to use
the Littlewood-Paley projections. Actually, we will arrive at a series of the form
$$\displaystyle\sum_{k=k_0}^{\infty}\sum_{j=j_0}^{\infty}\fd{||D^{\sigma}P_{j}
u(t)||_{{k}}^{{k}}}{2^{B_{k}}},\eqno{(1.6)}$$
where $\sigma$ is a fixed integer and $B_k$ is given by (5.1).

Our first step is to show that if (1.6) is convergent,
then it is bounded by a constant $\cal C$ depending only on $u_0$, $\nu$ and $T$. We call this
{\it uniform bound estimate}.
Now, to obtain the a priori estimate, we only need to show that (1.6) is always convergent.
Using the {\it uniform bound estimate}, we see (1.6) is always convergent on a closed time interval
and then it is left to show that if (1.6) is convergent at $T'$, then it is convergent on $[T',T'+\delta]$ with
some $\delta>0$. To do this, we separate (1.6) into low frequency part (finite $j$) and high frequency part (infinite $j$).

Our second step is to show the convergence of the low frequency part, that is,
$$\displaystyle\sum_{k=k_0}^{\infty}\sum_{j=j_0}^{J_0}\fd{||D^{\sigma}P_{j}
u(t)||_{{k}}^{{k}}}{2^{B_{k}}}\eqno{(1.7)}$$
is convergent on $[T',T'+\delta]$, where $J_0$ is a large integer. This is hard. We design a different series
$$\displaystyle\sum_{k=k_0}^{\infty}\sum_{j=j_0}^{J_0}\fd{||D^{\sigma}P_{j}
u(t)||_{{k}}^{{k}}}{2^{\hat B_{k}}}\eqno{(1.8)}$$
and show that (1.8) can not blow up before any given time $T$ if its the initial value is small enough
which can be satisfied by choosing $\hat B_k$ (defined by (6.2)) to be large enough. Then (1.8) is convergent
on $[T',T'+\delta]$ which implies the convergence of (1.7).
Devising suitable $B_k$ and $\hat B_k$ is the key to these two steps.

Our third step is to show a regularity improving result, from which it follows easily
the convergence of the high frequency part, that is, the convergence of
$$\displaystyle\sum_{k=k_0}^{\infty}\sum_{j=J_0+1}^{\infty}\fd{||D^{\sigma}P_{j}
u(t)||_{{k}}^{{k}}}{2^{B_{k}}}$$
on $[T',T'+\delta]$. This kind of regularity improving is not new essentially,
but we need a special form. These three steps are the scheme of our proof of the new a priori estimates.

We organize the paper as the following.
In Section 2, we study the Littlewood-Paley projections.
In Section 3, some well known results of
the Navier-Stokes equations are stated.
In Section 4, we give some preparations for our {\it attack}.
We demonstrate the above three steps of the proof of the new a priori estimates
in Section 5, 6 and 7 respectively. Then in section 8, we show our new a priori estimates
of the Navier-Stokes equations. The proof of Theorem 1.1 is given in the last section.

Throughout of this paper, the spacial dimension is confined to be 3
although all the main results can be extended to higher dimensions. We
will use standard notations in this paper.

$[x]$: the maximal integer less than or equal to $x$;

$||f||_p$: $L^p$ norm in $R^3$ of function $f$;

$||f(t)||_p$: $L^p$ norm in $R^3$ of function $f(x,t)$ defined on $R^3\times[0,T]$, if $t$ is clear from the context, it is simplified to be $||f||_p$;

$D^{\sigma}f$: all the $\sigma$-th order derivatives of $f$ with respect to space variables for any integer $\sigma\geq0$;

$P_j$: Littlewood-Paley projection;

$\otimes$: tensor product;

$C$: universal constants which may be different at different occurrence.

\bigskip

\section{The Littlewood-Paley projections}

In this section, we will study the Littlewood-Paley decomposition
which is an important tool for analysis (cf.[8] for more
discussion). For any test function $f$, define the Fourier
transformation and the inverse Fourier transformation by
$$\hat f(\xi)=\int_{R^3}f(x)e^{-2\pi ix\cdot\xi}dx$$ and
$$\check f(\xi)=\int_{R^3}f(x)e^{2\pi ix\cdot\xi}dx$$
respectively. Then we have $$\big(\hat f\big)^{\check{}}=\big(\check f\big)^{\hat{}}=f.$$

Let $\phi\in C^{\infty}(R^3)$ be a real radial function supported on
$B_2$ such that $\phi\equiv1$ on $B_1$ and $0\leq\phi\leq1$ in
$R^3$. Let $$\psi(\xi)=\phi(\xi)-\phi(2\xi)$$ and
$$\psi_j(\xi)=\psi(2^{-j}\xi)\mb{for}j=\dots,-3,-2,-1,0,1,2,3,\dots.$$
Then we have
$$\sum_{j=-\infty}^{+\infty}\psi_j(\xi)=1\mb{in}R^3\backslash\{0\}.$$

We now define the Littlewood-Paley projection $P_j$ by
$$\widehat{P_jf}(\xi)=\psi_j(\xi)\hat f(\xi).$$
Then we have
$$f=\sum_{j=-\infty}^{+\infty}P_jf.\eqno{(2.1)}$$
For simplicity, we denote
$$P_{\leq j}=\sum_{k=-\infty}^{j}P_k\mb{and}P_{\geq j}=\sum_{k=j}^{+\infty}P_k.$$

\bigskip

\lemma{2.1}{Let $1\leq q\leq+\infty$ and $-\infty<j<+\infty$. Then
$$||P_jf||_q\mb{and}||P_{\leq j}f||_q\leq C||f||_q,$$
where $C$ is a universal constant. }

\proof From
$$(P_jf)(x)=(\check\psi_j*f)(x)=\int_{R^3}\check\psi_j(y)f(x-y)dy$$
and Minkowski's inequality, we have
$$||P_jf||_q\leq\int_{R^3}|\check\psi_j((y)|||f||_qdy=||\check\psi_j||_1||f||_q.$$
Since $$||\check\psi_j||_1=||\check\psi_1||_1,$$ we see that
$$||P_jf||_q\leq||\check\psi_1||_1||f||_q.$$
Similarly, we have
$$||P_{\leq j}f||_q\leq||(\phi(2^{-j}\cdot))^{\check{}}||_1||f||_q=||\check{\phi}||_1||f||_q.$$
Let $$C=\max\{||\check\psi_1||_1,||\check{\phi}||_1\}$$ and the
proof of Lemma 2.1 is complete. $\Box$

\bigskip

From (2.1) and Lemma 2.1, we have the following so called cheap
Littlewood-Paley inequality (cf.[8]).

\bigskip

\theorem{2.2}{Let $1\leq q\leq+\infty$. Then
$$C\sup_{j}||P_jf||_q\leq||f||_q\leq\sum_{j=-\infty}^{+\infty}||P_jf||_q,$$
where $C$ is a universal constant.}

\bigskip

\lemma{2.3 (Bernstein's inequality)$^{[8]}$}{Let $1\leq q\leq
q'\leq+\infty$ and $j$ be an integer. Then
$$||P_jf||_{q'}\leq C2^{j(\frac3q-\frac3{q'})}||P_jf||_q,$$
where $C$ is a universal constant.}

\proof Let $\frac1{q'}=\frac{\lambda}{q}$ or $\lambda=\frac{q}{q'}$ and then
$$||P_jf||_{q'}\leq||P_jf||_{q}^{\lambda}||P_jf||_{\infty}^{1-\lambda}.\eqno{(2.2)}$$
From $\widehat{P_jf}=\psi_j\hat f=\phi(2^{-j-1}\cdot)\psi_j\hat f$, we deduce
$$P_jf=(\phi(2^{-j-1}\cdot))^{\check{}}*P_jf.$$
It follows that
$$||P_jf||_{\infty}\leq||(\phi(2^{-j-1}\cdot))^{\check{}}||_r||P_jf||_q=\left(2^{j+1}\right)^{3\frac{r-1}{r}}||\check\phi||_r||P_jf||_q,$$
where $\frac1r+\frac1q=1$ or $r=\frac{q}{q-1}$. It is easy to see
that $$\frac{r-1}{r}(1-\lambda)=\frac1q-\frac1{q'}.$$ In view of
(2.2),
$$||P_jf||_{q'}\leq||P_jf||_{q}^{\lambda}||P_jf||_{\infty}^{1-\lambda}\leq2^{(j+1)(\frac3q-\frac3{q'})}||\check\phi||_r^{1-\lambda}||P_jf||_q.$$
Since
$$\begin{array}{c}\displaystyle||\check\phi||_r=\left(\int_{R^3}\left|\check\phi(x)\frac{1+|x|^4}{1+|x|^4}\right|^rdx\right)^{\frac1r}
\leq\sup_{x\in R^3}\left(|\check\phi(x)|(1+|x|^4)\right)\left(\int_{R^3}\left(\frac{1}{1+|x|^4}\right)^rdx\right)^{\frac1r}\\[15pt]
\displaystyle\leq\sup_{x\in
R^3}\left(|\check\phi(x)|(1+|x|^4)\right)\left(\int_{R^3}\frac{1}{1+|x|^4}dx\right)^{\frac1r}\end{array}
$$
which can be bounded by a universal constant, we have the
conclusion. $\Box$

\bigskip

\lemma{2.4}{Let $j$ be an integer, $f$ and $g$ be two test functions. We have the following product inequality,
$$\begin{array}{c}\displaystyle |P_j(fg)|\displaystyle
\leq
\\[15pt]\displaystyle\mbox{}\hskip0cm
\left|P_j\Bigg\{\left(\sum_{m=-\infty}^{j-3}P_mf\right)\left(\sum_{m'=j-2}^{j+2}P_{m'}g\right)\Bigg\}\right|
+
\left|P_j\Bigg\{\left(\sum_{m=j-2}^{j+2}P_mf\right)\left(\sum_{m'=-\infty}^{j-3}P_{m'}g\right)\Bigg\}\right|
\\[15pt]\displaystyle\mbox{}\hskip0cm
+\left|P_j\Bigg\{\left(\sum_{m=j-2}^{j+2}P_mf\right)
\left(\sum_{m'=j-2}^{j+2}P_{m'}g\right)\Bigg\}\right|+
\\[15pt]\displaystyle\mbox{}\hskip0cm
\left|P_j\Bigg\{\left(\sum_{m=j+3}^{\infty}P_mf\right)
\left(\sum_{m'=m-3}^{m+3}P_{m'}g\right)\Bigg\}\right|
+\left|P_j\Bigg\{\left(\sum_{m'=j+3}^{\infty}P_{m'}g\right)\left(\sum_{m=m'-3}^{m'+3}P_mf\right)\Bigg\}\right|.
\end{array}$$
}

\proof
From (2.1), it follows that
$$\begin{array}{l}\displaystyle P_j(fg)\displaystyle=P_j\left(\left(\sum_{m=-\infty}^{\infty}P_mf\right)
\left(\sum_{m'=-\infty}^{\infty}P_{m'}g\right)\right)
\\[15pt]\displaystyle\mbox{}\hskip1cm
=P_j\Bigg\{\left(\sum_{m=-\infty}^{j-3}P_mf+\sum_{m=j-2}^{j+2}P_mf+\sum_{m=j+3}^{\infty}P_mf\right)
\\[15pt]\displaystyle\displaystyle\mbox{}\hskip2.2cm\times
\left(\sum_{m'=-\infty}^{j-3}P_{m'}g+\sum_{m'=j-2}^{j+2}P_{m'}g+\sum_{m'=j+3}^{\infty}P_{m'}g\right)\Bigg\}.
\end{array}$$
Since
$$P_j\left(P_{\leq j-3}fP_{\leq j-3}g\right)=0$$
and
$$P_j\left(P_mfP_{m'}g\right)=0$$
as $m\geq j+3$ and $|m-m'|>3$, we have the conclusion clearly.
$\Box$

\bigskip

\lemma{2.5}{Let $1=j_0\leq j$ be two integers,  $1\leq q\leq\infty,2\leq q_0, q_1\leq\infty$ be
three real numbers, and $f$ and $g$ be two test functions.
If $\frac1q=\frac1{q_0}+\frac1{q_1}$, then
$$\begin{array}{l}\displaystyle||P_j(fg)||_q\leq\displaystyle C\Bigg\{
\alpha_j||f||_2||g||_2
+
\sum_{m=\max\{j_0,j-2\}}^{\infty}||P_mf||_{q_1}||g||_{q_0}
\\[15pt]\displaystyle\mbox{}\hskip4cm
+
\sum_{m'=\max\{j_0,j-2\}}^{\infty}||P_{m'}g||_{q_1}||f||_{q_0}
\Bigg\},\end{array}$$ where $C$ is a universal constant and
$$\alpha_j=\left\{\begin{array}{l}1\mb{as}j=1,2;\\[15pt]
0\mb{as}j\geq3.\end{array}\right.\eqno{(2.3)}$$

. }

\proof From Lemma 2.4 and H\"{o}lder's inequality, it follows that
$$\begin{array}{c}\displaystyle ||P_j(fg)||_{q}\displaystyle
\leq ||P_{\leq j-3}f||_{q_0}\sum_{m'=j-2}^{j+2}||P_{m'}g||_{q_1}
+||P_{\leq j-3}g||_{q_0}\sum_{m=j-2}^{j+2}||P_{m}f||_{q_1}
\\[15pt]\displaystyle\mbox{}\hskip0cm
+\sum_{m=j-2}^{j+2}\sum_{m'=j-2}^{j+2}||P_mf||_{q_0}||P_{m'}g||_{q_1}+
\\[15pt]\displaystyle\mbox{}\hskip0cm
\sum_{m=j+3}^{\infty}\sum_{m'=m-3}^{m+3}||P_mf||_{q_1}||P_{m'}g||_{q_0}
+
\sum_{m'=j+3}^{\infty}\sum_{m=m'-3}^{m'+3}||P_mf||_{q_0}||P_{m'}g||_{q_1}
.
\end{array}\eqno{(2.4)}$$
From Lemma 2.1 and 2.3, we have
$$\begin{array}{c}\displaystyle
||P_{\leq j-3}f||_{q_0}\sum_{m'=j-2}^{j+2}||P_{m'}g||_{q_1}=
||P_{\leq j-3}f||_{q_0}\left(\sum_{m'=j-2}^{\max\{j_0,j-2\}-1}+\sum_{m'=\max\{j_0,j-2\}}^{j+2}\right)||P_{m'}g||_{q_1}
\\[15pt]\displaystyle\mbox{}\hskip0cm
\leq C\left(\alpha_j||f||_2||g||_2+
||f||_{q_0}\sum_{m'=\max\{j_0,j-2\}}^{j+2}||P_{m'}g||_{q_1}\right),
\end{array}$$
where $C$ is a universal constant and $\alpha_j$ is given by (2.3).
Similarly,
$$\begin{array}{c}\displaystyle
||P_{\leq j-3}g||_{q_0}\sum_{m=j-2}^{j+2}||P_{m}f||_{q_1}
\leq C\left(\alpha_j||f||_2||g||_2+
||g||_{q_0}\sum_{m=\max\{j_0,j-2\}}^{j+2}||P_{m}f||_{q_1}\right)
\end{array}$$
and
$$\begin{array}{c}\displaystyle
\sum_{m=j-2}^{j+2}\sum_{m'=j-2}^{j+2}||P_mf||_{q_0}||P_{m'}g||_{q_1}
\leq C\left(\alpha_j||f||_2||g||_2+
||f||_{q_0}\sum_{m'=\max\{j_0,j-2\}}^{j+2}||P_{m'}g||_{q_1}\right)
\end{array}$$
where $C$ are universal constants and $\alpha_j$ is given by (2.3).
From Lemma 2.1, we have
$$\begin{array}{c}\displaystyle
\sum_{m=j+3}^{\infty}\sum_{m'=m-3}^{m+3}||P_mf||_{q_1}||P_{m'}g||_{q_0}
+
\sum_{m'=j+3}^{\infty}\sum_{m=m'-3}^{m'+3}||P_mf||_{q_0}||P_{m'}g||_{q_1}
\\[15pt]\displaystyle\mbox{}\hskip0cm
\leq C\left(
\sum_{m=j+3}^{\infty}||P_mf||_{q_1}||g||_{q_0}
+
\sum_{m'=j+3}^{\infty}||f||_{q_0}||P_{m'}g||_{q_1}\right),
\end{array}$$
where $C$ is a universal constant.
Plug the above inequalities into (2.4) and then we have the conclusion.
$\Box$

\bigskip

The following lemma gives one of the key reason why the Littlewood-Paley projection is useful.
We refer to [8] for its proof.

\bigskip

\lemma{2.6}{Let $1\leq q\leq\infty$ be a real number and $j$ be an integer. Then we
have
$$C_12^{j}||P_jf||_q\leq||\nabla P_jf||_q\leq C_22^{j}||P_jf||_q,$$
where $C_1$ and $C_2$ are universal constants.}

\bigskip

\lemma{2.7}{Let $1=j_0\leq j$, $0\leq\sigma$ be three integers, $1\leq q\leq\infty,2\leq q_0, q_1\leq\infty$ be
three real numbers, and $f$ and $g$ be two test functions.
If $\frac1q=\frac1{q_0}+\frac1{q_1}$, then
$$\begin{array}{l}\displaystyle\mbox{}\hskip0cm|| D^{\sigma}P_j(fg)||_q\leq
C(\sigma)\Bigg\{\alpha_j||f||_2||g||_{2}+
\sum_{m=\max\{j_0,j-2\}}^{+\infty}2^{(j-m)\sigma}|| D^{\sigma}P_mf||_{q_1}||g||_{q_0}
\\[15pt]\displaystyle\mbox{}\hskip4cm
+
\sum_{m'=\max\{j_0,j-2\}}^{+\infty}2^{(j-m')\sigma}|| D^{\sigma}P_{m'}g||_{q_1}||f||_{q_0}
\Bigg\}
,\end{array}$$ where
$C(\sigma)$ is a constant depending only on $\sigma$ and $\alpha_j$ is given by (2.3).}

\proof From Lemma 2.6, we have
$$|| D^{\sigma}P_j(fg)||_q\leq C(\sigma)2^{j\sigma}||P_j(fg)||_q,\eqno{(2.5)}$$
where $C(\sigma)$ is a constant depending only on $\sigma$.
From Lemma 2.5 and 2.6, we have
$$\begin{array}{l}\displaystyle||P_j(fg)||_q\leq\displaystyle C\Bigg\{
\alpha_j||f||_2||g||_2
+
\sum_{m=\max\{j_0,j-2\}}^{\infty}||P_mf||_{q_1}||g||_{q_0}
\\[15pt]\displaystyle\mbox{}\hskip4cm
+
\sum_{m'=\max\{j_0,j-2\}}^{\infty}||P_{m'}g||_{q_1}||f||_{q_0}
\Bigg\}
\\[15pt]\displaystyle\mbox{}\hskip0cm
\leq\displaystyle C(\sigma)\Bigg\{
\alpha_j||f||_2||g||_2
+
\sum_{m=\max\{j_0,j-2\}}^{\infty}2^{-m\sigma}||D^{\sigma}P_mf||_{q_1}||g||_{q_0}
\\[15pt]\displaystyle\mbox{}\hskip4cm
+
\sum_{m'=\max\{j_0,j-2\}}^{\infty}2^{-m'\sigma}||D^{\sigma}P_{m'}g||_{q_1}||f||_{q_0}
\Bigg\},\end{array}$$
where $C$ is a universal constant,
$C(\sigma)$ is a constant depending only on $\sigma$, and $\alpha_j$ is given by (2.3).
Combining it with (2.5), we see the conclusion clearly.
$\Box$

\bigskip

\corollary{2.8}{Suppose all the assumptions of Lemma 2.7 hold. Let $\hat q>1$ be a real number. Then
$$\begin{array}{l}\displaystyle\mbox{}\hskip0cm|| D^{\sigma}P_j(fg)||_q^{\hat q}\leq
C(\sigma)^{\hat q}\Bigg\{\alpha_j||f||_2^{\hat q}||g||_{2}^{\hat q}+
\sum_{m=\max\{j_0,j-2\}}^{\infty}2^{(j-m)\sigma}|| D^{\sigma}P_mf||_{q_1}^{\hat q}||g||_{q_0}^{\hat q}
\\[15pt]\displaystyle\mbox{}\hskip4cm
+
\sum_{m'=\max\{j_0,j-2\}}^{\infty}2^{(j-m')\sigma}|| D^{\sigma}P_{m'}g||_{q_1}^{\hat q}||f||_{q_0}^{\hat q}
\Bigg\}
,\end{array}$$ where
$C(\sigma)$ is a constant depending only on $\sigma$ and $\alpha_j$ is given by (2.3).}

\proof From Lemma 2.7, we have
$$\begin{array}{l}\displaystyle\mbox{}\hskip0cm|| D^{\sigma}P_j(fg)||_q^{\hat q}\leq
C(\sigma)^{\hat q}\Bigg\{\alpha_j||f||_2||g||_{2}+
\sum_{m=\max\{j_0,j-2\}}^{+\infty}2^{(j-m)\sigma}|| D^{\sigma}P_mf||_{q_1}||g||_{q_0}
\\[15pt]\displaystyle\mbox{}\hskip4cm
+
\sum_{m'=\max\{j_0,j-2\}}^{+\infty}2^{(j-m')\sigma}|| D^{\sigma}P_{m'}g||_{q_1}||f||_{q_0}
\Bigg\}^{\hat q}
,\end{array}\eqno{(2.6)}$$ where
$C(\sigma)$ is a constant depending only on $\sigma$ and $\alpha_j$ is given by (2.3).
From H\"{o}lder's inequality, we have
$$\begin{array}{l}\displaystyle\mbox{}\hskip0cm\alpha_j||f||_2||g||_{2}+
\sum_{m=\max\{j_0,j-2\}}^{+\infty}2^{(j-m)\sigma}|| D^{\sigma}P_mf||_{q_1}||g||_{q_0}
\\[15pt]\displaystyle\mbox{}\hskip3cm
+
\sum_{m'=\max\{j_0,j-2\}}^{+\infty}2^{(j-m')\sigma}|| D^{\sigma}P_{m'}g||_{q_1}||f||_{q_0}
\\[15pt]\displaystyle\mbox{}\hskip0cm
\leq\Bigg\{\alpha_j+\sum_{m=\max\{j_0,j-2\}}^{+\infty}2^{(j-m)\sigma}+\sum_{m=\max\{j_0,j-2\}}^{+\infty}2^{(j-m)\sigma}
\Bigg\}^{\frac1{\hat q'}}
\\[15pt]\displaystyle\mbox{}\hskip1cm\times
\Bigg\{\alpha_j||f||_2^{\hat q}||g||_{2}^{\hat q}+
\sum_{m=\max\{j_0,j-2\}}^{+\infty}2^{(j-m)\sigma}|| D^{\sigma}P_mf||_{q_1}^{\hat q}||g||_{q_0}^{\hat q}
\\[15pt]\displaystyle\mbox{}\hskip3cm
+
\sum_{m'=\max\{j_0,j-2\}}^{+\infty}2^{(j-m')\sigma}|| D^{\sigma}P_{m'}g||_{q_1}^{\hat q}||f||_{q_0}^{\hat q}\Bigg\}^{\frac1{\hat q}}
\\[15pt]\displaystyle\mbox{}\hskip0cm
\leq\left(1+ 8^{\sigma}\right)
\Bigg\{\alpha_j||f||_2^{\hat q}||g||_{2}^{\hat q}+
\sum_{m=\max\{j_0,j-2\}}^{+\infty}2^{(j-m)\sigma}|| D^{\sigma}P_mf||_{q_1}^{\hat q}||g||_{q_0}^{\hat q}
\\[15pt]\displaystyle\mbox{}\hskip3.6cm
+
\sum_{m'=\max\{j_0,j-2\}}^{+\infty}2^{(j-m')\sigma}|| D^{\sigma}P_{m'}g||_{q_1}^{\hat q}||f||_{q_0}^{\hat q}\Bigg\}^{\frac1{\hat q}}
,\end{array}$$
where $\frac1{\hat q}+\frac1{\hat q'}=1$. Combining it with (2.6), we have the conclusion.
$\Box$

\bigskip

\section{Classical results}

In this section, we state some classical results of
Navier-Stokes equations, which are all well known.
Although we can not find exact references for some of them, we still omit the proofs here.
We refer to [9] for the definition of Leray's weak solutions and strong solutions of Navier-Stokes equations.
Here the domain of the space is $R^3$. We should point out that any strong solutions will be Leray's weak solutions.

\bigskip

\theorem{3.1}{Let $u$ be a Leray's weak solution of
(1.1) with $f=0$ and $u(0)=u_0$
($\mbox{div}u_0=0$). Then
$$\displaystyle||u(t)||_{2}\leq||u_0||_{2}\mb{and}\int_0^{t}||\nabla
u(s)||_{2}^2ds\leq\frac1{2\nu}||u_0||_{2}^2$$ for any $t\geq0$.}

\bigskip

\corollary{3.2}{Let $u$ be given by Theorem 3.1. Then
$$\int_0^{t}||u(s)||_4^{\frac83}ds\leq
\frac{C}{\nu}||u_0||_{2}^{\frac83}$$
for any $t\geq0$, where $C$ is a universal constant. }

\proof From interpolation
inequality and Sobolev's embedding inequality,
$$||u||_4\leq ||u||_2^{\frac14}||u||_6^{\frac34}\leq C||u||_2^{\frac14}||\nabla u||_2^{\frac34},$$
where $C$ is a universal constant.
Therefore
$$\int_0^{t}||u(s)||_4^{\frac83}ds\leq
C\int_0^{t}||u(s)||_2^{\frac23}||\nabla u(s)||_2^{2}ds.$$
From Theorem 3.1,
we deduce the conclusion easily. $\Box$

\bigskip

\theorem{3.3 (Uniqueness)}{Let $u_0\in H^1(R^3)$ be a divergence-free vector field.
Then the strong solution of (1.1) with $u(0)=u_0$ and $f\equiv0$ is unique. }

\bigskip

\theorem{3.4 (Short time existence)}{Let $u_0\in H^1(R^3)$ be a divergence-free vector field.
Then there exists $T^*>0$ depending only on $\nu$ and $||\nabla u_0||_2$
such that the strong solution of (1.1) exists on the time interval $[0,T^*]$ with $u(0)=u_0$ and $f\equiv0$.
 }

\bigskip

\theorem{3.5 (Regularity)}{Suppose $u_0$ is a smooth, divergence-free vector field
in $R^3$ satisfying (1.2). Let
$u$ and $p$ be the strong solution of (1.1) on the time interval $[0,T]$ with $u(0)=u_0$ and $f\equiv0$.
Then $u$ and $p$ are smooth and
$u,p\in C^1\left([0,T],W^{m,2}(R^3)\cap W^{m,\infty}(R^3)\right)$ for any integer $m\geq1$.
}

\bigskip

\theorem{3.6 (Blow up)}{Let $u_0\in H^1(R^3)$ be a divergence-free vector field.
Suppose $0<T<\infty$ and $[0,T)$ is the largest time interval that (1.1) has
strong solution with $u(0)=u_0$ and $f\equiv0$. Then
$$\limsup_{t\rightarrow T^-}||u(t)||_{\infty}=\infty.$$

}

\bigskip

\remark{3.7} In Theorem 3.5, if we do not assume that $u_0$ is smooth and satisfies (1.2), then we will have
for any $0<\hat T<T$, $u$ and $p$ are smooth on the time interval $[\hat T,T]$ and
$u,p\in C^1\left([\hat T,T],W^{m,2}(R^3)\cap W^{m,\infty}(R^3)\right)$ for any integer $m\geq1$. $\Box$

\bigskip

Finally, we prove a simple lemma.

\bigskip

\lemma {3.8}{Let $q\geq2$ be a real number and
$\sigma$ and $j$ be two positive integers. Suppose $u$ and $p$ satisfy (1.1) with $f\equiv0$.  Then
$$||D^{\sigma}P_jp||_{q}\leq Cq||D^{\sigma}P_j(u\otimes u)||_{q},\eqno{(3.1)}$$
where $C$ is a universal constant.
}

\proof Take divergence on the both sides of the first equation of (1.1) and then
$$\sum_{i,j=1}^{3}\partial_i\partial_j(u_iu_j)+\Delta p=0.$$
Now take the operator $D^{\sigma}P_j$ on the both sides and then
$$\sum_{i,j=1}^{3}D^{\sigma}P_j\partial_i\partial_j(u_iu_j)+\Delta D^{\sigma}P_j p=0.$$
From Calder\'{o}n-Zygmund's estimate, we have (3.1).
$\Box$

\bigskip

\section{Some lemmas}

In this section, we will show some lemmas which will be used to prove our new a priori estimates.

\bigskip

\lemma{4.1}{Let $k\geq2$ and suppose $f\in W^{2,2}(R^3)\cap W^{2,\infty}(R^3)$. Then we have
$$||f||_{3k}^k\leq Ck^2\int_{R^3}|f|^{k-2}|\nabla f|^2,$$
where $C$ is a universal constant.
}

\proof From Sobolev's embedding inequality, we have
$$\begin{array}{l}\displaystyle
||f||_{3k}^k=|||f|^{\frac{k}2}||_{6}^2\leq C||\nabla|f|^{\frac{k}2}||_2^2
\leq Ck^2\int_{R^3}|f|^{k-2}|\nabla f|^2,
\end{array}$$
where $C$ is a universal constant.
$\Box$

\bigskip

\lemma{4.2}{Let $T>0$ be a real number and $\{f_k(t)\}_{k=1}^{\infty}$ and $\{g_k(t)\}_{k=1}^{\infty}\subset C^1[0,T]$ be two nonnegative function sequences. Suppose
$$\sum_{k=1}^{\infty}f_k(t)\leq{\cal B}\mb{and}
\frac{d}{dt}f_k(t)\leq g_k(t),\forall0\leq t\leq T\mb{and}k\geq1,\eqno{(4.1)}$$
where ${\cal B}$ is a constant. Then we have
$$\frac{d}{dt}\sum_{k=1}^{\infty}f_k(t)\leq\sum_{k=1}^{\infty}g_k(t)$$
for any $0\leq t\leq T$.
}

\proof From  Lebesgue's dominated convergence theorem, (4.1) and Fatou's Lemma,
we deduce
$$
\begin{array}{l}\displaystyle
\int_0^T\phi(t)\frac{d}{dt}\sum_{k=1}^{\infty}f_k(t)
=-\int_0^T\sum_{k=1}^{\infty}f_k(t)\frac{d}{dt}\phi(t)
=-\sum_{k=1}^{\infty}\int_0^Tf_k(t)\frac{d}{dt}\phi(t)\\[15pt]\displaystyle
=\sum_{k=1}^{\infty}\int_0^T\phi(t)\frac{d}{dt}f_k(t)
\leq\sum_{k=1}^{\infty}\int_0^T\phi(t)g_k(t)
\leq\int_0^T\phi(t)\sum_{k=1}^{\infty}g_k(t)
\end{array}
$$
for any test function $\phi\geq0$. This implies the conclusion clearly.
$\Box$

\bigskip

\lemma{4.3}{Let $\epsilon,T,{\cal B}>0$ and $M>1$ be real numbers.
Let $0\leq F\in C^1[0,T]$ and $0\leq g\in C[0,T]$. Suppose
$$\int_0^Tg(t)\leq{\cal B},\quad F(0)\leq\epsilon$$ and
$$\frac{d}{dt}F(t)\leq g(t)\left(\epsilon+F(t)+F^M(t)\right).\eqno{(4.2)}$$
If $$\epsilon\leq\frac{1}{\left(3e^{\cal B}\right)^{\frac{M}{M-1}}},\eqno{(4.3)}$$
then we have $$F(t)\leq3\epsilon e^{\cal B}\eqno{(4.4)}$$
for any $0\leq t\leq T$.
}

\proof Suppose by the contradiction that (4.4) is not true. Then since
$$F(0)\leq\epsilon<3\epsilon e^{\cal B},$$
there exists $T'\in(0,T)$ such that (4.4) holds for $t\in[0,T']$ and
$$F(T')=3\epsilon e^{\cal B}.\eqno{(4.5)}$$

It is easy to see that (4.3) and (4.4) which we assume holds for $t\in[0,T']$ imply that
$$F^M(t)\leq\epsilon\mb{for}t\in[0,T'].$$
From (4.2), it follows that
$$\frac{d}{dt}F(t)\leq g(t)\left(2\epsilon+F(t)\right)\mb{for}t\in[0,T']$$
or
$$\ln\frac{2\epsilon+F(t)}{2\epsilon+F(0)}\leq\int_0^{T'}g(t)\leq{\cal B}\mb{for}t\in[0,T'].$$
Therefore
$$F(T')\leq\left(2\epsilon+F(0)\right)e^{\cal B}-2\epsilon<3\epsilon e^{\cal B}.$$
This contradicts with (4.5).
$\Box$

\bigskip

\lemma{4.4\footnote{This lemma is given by Prof.Lihe Wang.}}{Let $k\geq2$ and $\sigma\geq1$ be integers.
Suppose $f\in W^{2,2}(R^3)\cap W^{2,\infty}(R^3)$ . Then we have
$$\int_{R^3}|D^{\sigma} f|^k\leq C(\sigma)k^2\left(\int_{R^3}|D^{\sigma} f|^{k-2}|D^{\sigma+1}f|^2\right)^{\frac{k}{k+2}}
\left(\int_{R^3}|D^{\sigma-1}f|^{k}\right)^{\frac{2}{k+2}},\eqno{(4.6)}$$
where $C(\sigma)$ is a constant depending only on $\sigma$.
}

\proof We only need to prove (4.6) with the assumption that $f\in C_c^{\infty}(R^3)$.
From Green's and H\"{o}lder's formula, we see
$$\begin{array}{c}\displaystyle
\int_{R^3}|D^{\sigma} f|^k
=-\int_{R^3}D^{\sigma-1}fD\left(|D^{\sigma} f|^{k-2}D^{\sigma} f\right)
\\[15pt]\displaystyle\mbox{}\hskip0.5cm
\leq C(\sigma)
(k-1)\int_{R^3}|D^{\sigma-1}f||D^{\sigma} f|^{k-2}|D^{\sigma+1} f|
\\[15pt]\displaystyle\mbox{}\hskip0.5cm
\leq C(\sigma)(k-1)\left(\int_{R^3}|D^{\sigma-1}f|^k\right)^{\frac1k}
\left(\int_{R^3}|D^{\sigma} f|^{k}\right)^{\frac{k-2}{2k}}
\left(\int_{R^3}|D^{\sigma} f|^{k-2}|D^{\sigma+1} f|^2\right)^{\frac12}
,\end{array}$$ where $C(\sigma)$ is a constant depending only on $\sigma$.
It follows (4.6) clearly.
$\Box$

\bigskip

\lemma{4.5}{Let $k\geq2$, $\sigma\geq1$ and $j$ be integers.
Suppose $f\in W^{2,2}(R^3)\cap W^{2,\infty}(R^3)$. Then we have
$$\int_{R^3}|D^{\sigma} P_j f|^k\leq C(\sigma)k^22^{-2j}\int_{R^3}|D^{\sigma} P_j f|^{k-2}|D^{\sigma+1} P_jf|^2
,\eqno{(4.7)}$$
where $C(\sigma)$ is a constant depending only on $\sigma$.
}

\proof From Lemma 4.4, we have
$$\int_{R^3}|D^{\sigma} P_j f|^k\leq C(\sigma)k^2\left(\int_{R^3}|D^{\sigma} P_j f|^{k-2}|D^{\sigma+1}P_jf|^2\right)^{\frac{k}{k+2}}
\left(\int_{R^3}|D^{\sigma-1}P_jf|^{k}\right)^{\frac{2}{k+2}}.$$
From Lemma 2.6, we have
$$\left(\int_{R^3}|D^{\sigma-1}P_jf|^{k}\right)^{\frac{2}{k+2}}\leq C2^{-\frac{2k}{k+2}j}\left(\int_{R^3}|D^{\sigma} P_jf|^{k}\right)^{\frac{2}{k+2}},$$
where $C$ is a universal constant.
Then we see (4.7) clearly.
$\Box$

\bigskip

\lemma{4.6}{
Let $k_0\geq1$, $\sigma\geq0$ be two integers, $T>0$ be a real number and $u\in C\left([0,T],W^{\sigma+1,k}(R^3)\right)$ for any $k\geq k_0$. Then for any $k\geq k_0$,
$$\displaystyle\sum_{j=j_0}^{\infty}||D^{\sigma}P_{j}u(t)||_{k}^{k}$$ is continuous as a function on $[0,T]$,
and there exists a constant ${\cal B}>0$ such that
$$\displaystyle\sum_{j=j_0}^{\infty}||D^{\sigma}P_{j}u||_{k}^{k}\leq{\cal B}$$
for any $0\leq t\leq T$.
}

\proof From Lemma 2.6, we have
$$\begin{array}{l}\displaystyle
||D^{\sigma}P_{j}u(t)||_{k}\leq C2^{-j}||D^{\sigma+1}P_{j}u(t)||_{k}
\leq  C2^{-j}\sup_{t\in[0,T]}||D^{\sigma+1}u(t)||_{k}.
\end{array}$$
Since
$$\begin{array}{l}\displaystyle
\sum_{j=j_0}^{\infty}\left(C2^{-j}\sup_{t\in[0,T]}||D^{\sigma+1}u(t)||_{k}\right)^k
\end{array}$$
is convergent, we see the conclusion clearly.
$\Box$

\bigskip

\lemma{4.7}{Let $k_0\geq10$, $j_0\geq1$ and $\sigma\geq1$ be three integers, $B>0$ be a real number and $u_0$ be a function satisfying (1.2).
Then there exists $\tilde B_0>0$ depending only on $u_0$ and $\sigma$ such that if $B\geq \tilde B_0$, then
$$\sum_{k=k_0}^{\infty}\sum_{j=j_0}^{\infty}\frac{||D^{\sigma}P_j
u_0||_{k}^{k}}{2^{B(1-\frac{1}{\sqrt{k}})k}}\leq2^{-\frac{B}4}.$$
}

\proof From Lemma 2.6, we have
$$\begin{array}{l}\displaystyle
\sum_{k=k_0}^{\infty}\sum_{j=j_0}^{\infty}\frac{||D^{\sigma}P_j
u_0||_{k}^{k}}{2^{B(1-\frac{1}{\sqrt{k}})k}}\leq
\sum_{k=k_0}^{\infty}\sum_{j=j_0}^{\infty}\frac{\left(C2^{-j}||D^{\sigma+1}P_j
u_0||_{k}\right)^{k}}{2^{B(1-\frac{1}{\sqrt{k}})k}}
\\[15pt]\displaystyle\mbox{}\hskip0.1cm
\leq
\sum_{k=k_0}^{\infty}\sum_{j=j_0}^{\infty}2^{-jk}\frac{\left(C||D^{\sigma+1}u_0||_{k}\right)^{k}}{2^{B(1-\frac{1}{\sqrt{k}})k}}
\leq
\sum_{k=k_0}^{\infty}\frac{\left(C||D^{\sigma+1}u_0||_{k}\right)^{k}}{2^{B(1-\frac{1}{\sqrt{k}})k}},
\end{array}$$
where $C$ is universal constant. Let
$$\tilde B_0=2C\left(||D^{\sigma+1}u_0||_{2}+||D^{\sigma+1}u_0||_{\infty}\right)+4.$$
By interpolation inequality,
$$\begin{array}{l}\displaystyle
||D^{\sigma+1}u_0||_{k}\leq||D^{\sigma+1}u_0||_{2}^{\frac2k}||D^{\sigma+1}u_0||_{\infty}^{1-\frac2k}
\leq
||D^{\sigma+1}u_0||_{2}+||D^{\sigma+1}u_0||_{\infty}
\\[15pt]\displaystyle\mbox{}\hskip2cm
\leq\frac{\tilde B_0-4}{2C}
,
\end{array}$$
we see that if $B\geq \tilde B_0$,
$$\begin{array}{l}\displaystyle
\sum_{k=k_0}^{\infty}\frac{\left(C||D^{\sigma+1}u_0||_{k}\right)^{k}}{2^{B(1-\frac{1}{\sqrt{k}})k}}
=
\sum_{k=k_0}^{\infty}\frac1{2^{\frac14Bk}}
\left(\frac{C||D^{\sigma+1}u_0||_{k}}{2^{\frac34B-\frac{1}{\sqrt{k}}}}\right)^k
\\[15pt]\displaystyle\mbox{}\hskip2cm
\leq
\sum_{k=k_0}^{\infty}\frac1{2^{\frac14Bk}}
\left(\frac{B/2}{2^{B/2}}\right)^k
\leq2^{-\frac{B}4},
\end{array}$$
where $B\geq\tilde B_0\geq4$ is used.
$\Box$

\bigskip

\section{Uniform bound estimate}

In this section, we will prove the following {\it uniform bound estimate}, Theorem 5.1, which is the first key step to show our new
a priori estimates. We design the following series
$$\displaystyle\sum_{k=k_0}^{\infty}\sum_{j=j_0}^{\infty}\fd{||D^{\sigma}P_{j}
u(t)||_{{k}}^{{k}}}{2^{B_{k}}},\eqno{(5.1)}$$
where $j_0=1$, $k_0=100$ and the power $$B_k=\left(B+1+\frac{1}{\sqrt{k}}\right)k\eqno{(5.2)}$$
with positive constant $B$ to be given later.

For convenience, we define the following Condition (S) of $u$ and $p$:
$$\left\{\begin{array}{l}\displaystyle
(i)\hsa u\mb{and}p\in C^1\left([0,T],W^{m,2}(R^3)\cap W^{m,\infty}(R^3)\right)
\\[8pt]\mbox{}\hskip0.3cm\mb{for any integer}m\geq1;\\[8pt]
(ii)\hsa u\mb{and}p\mb{satisfy (1.1) with}f\equiv0\mb{and}u(x,0)=u_0(x),
\end{array}\right.\leqno{(S)}$$
for the given real number $T>0$ and the given suitable function $u_0$.

Our main theorem of this section is:

\bigskip

\theorem{5.1 (Uniform bound)}{
Let $\sigma=2$, $j_0=1$, $k_0=100$, $0<T'\leq T$ and $B>0$ be real numbers, $B_k$ be given by (5.1) for any $k\geq k_0$
and $u_0$ be a function satisfying (1.2).
Suppose $u$ and $p$ satisfy Condition (S). If
$$\displaystyle\sum_{k=k_0}^{\infty}\sum_{j=j_0}^{\infty}\fd{||D^{\sigma}P_{j}
u(t)||_{k}^{k}}{2^{B_k}}\leq{\cal B}\eqno{(5.3)}$$
for any $0\leq t\leq T'$, where ${\cal B}$ is a constant, then we have
$$\begin{array}{l}\displaystyle
\sum_{k=k_0}^{\infty}\sum_{j=j_0}^{\infty}
\frac{||D^{\sigma}P_j
u(t)||_{k}^{k}}{2^{B_k}}
\leq
{\cal C}
\sum_{k=k_0}^{\infty}\sum_{j=j_0}^{\infty}\frac{||D^{\sigma}P_j
u_0||_{k}^{k}}{2^{B_k}}+{\cal C}-1
\end{array}\eqno{(5.4)}$$
with
$$\displaystyle{\cal C}=e^{C\left(1+\nu^{-2}\right)\left((1+||u_0||_2^5)T+\nu^{-1}||u_0||_2^{\frac83}\right)}\eqno{(5.5)}$$
for any $0\leq t\leq T'$, where
$C>0$ is a universal constant.}

\bigskip

\proof Let $j\geq j_0$, $k\geq k_0$ and $0\leq t\leq T'$. We divide the proof into 6 steps.

\bigskip

{\it  Step 1.}\hsa
We first take the Littlewood-Paley projection on both sides
of Navier-Stokes equation:
$$P_ju_{t}-\nu P_j\Delta u+P_j((u\cdot \nabla)u)+P_j\nabla p=0$$
and then use the differential operator $D^{\sigma}$ on both
sides:
$$D^{\sigma}P_ju_{t}-\nu D^{\sigma}P_j\Delta u+D^{\sigma}P_j((u\cdot \nabla)u)+D^{\sigma}P_j\nabla p=0.$$
Finally, multiply $|D^{\sigma}P_ju|^{k-2}D^{\sigma}P_ju$ and then
integrate over $R^3$ on both sides:
$$\begin{array}{l}\displaystyle\mbox{}\hskip-0.5cm\int_{R^3}|D^{\sigma}P_ju|^{k-2}D^{\sigma}P_juD^{\sigma}P_ju_{t}
\\[15pt]\displaystyle
-\nu\int_{R^3}|D^{\sigma}P_ju|^{k-2}D^{\sigma}P_juD^{\sigma}P_j\Delta
u\\[15pt]
\mbox{}\hskip0.5cm+\displaystyle\int_{R^3}|D^{\sigma}P_ju|^{k-2}D^{\sigma}P_juD^{\sigma}P_j((u\cdot
\nabla)u) \\[15pt]
\displaystyle
\mbox{}\hskip1.5cm+\int_{R^3}|D^{\sigma}P_ju|^{k-2}D^{\sigma}P_juD^{\sigma}P_j\nabla
p=0.\end{array}\eqno{(5.6)}$$

It is easy to see that
$$\int_{R^3}|D^{\sigma}P_ju|^{k-2}D^{\sigma}P_juD^{\sigma}P_ju_{t}=
\frac1{k}\frac{d}{dt}\int_{R^3}|D^{\sigma}P_ju|^{k}.\eqno{(5.7)}$$

From Condition (S), by Green's formula,
$$\begin{array}{l}\displaystyle\mbox{}\hskip-0.5cm-\int_{R^3}|D^{\sigma}P_ju|^{k-2}D^{\sigma}P_juD^{\sigma}P_j\Delta u
=\int_{R^3}|D^{\sigma}P_ju|^{k-2}|D^{\sigma}P_j\nabla u|^2
\\[15pt]\displaystyle\mbox{}\hskip1cm
+
(k-2)\int_{R^3}|D^{\sigma}P_ju|^{k-2}|\nabla|D^{\sigma}P_ju||^2.\end{array}\eqno{(5.8)}$$
Similarly, we also have
$$\begin{array}{l}\displaystyle\mbox{}\hskip-0.5cm\int_{R^3}|D^{\sigma}P_ju|^{k-2}D^{\sigma}P_juD^{\sigma}P_j(u\cdot \nabla)u=
\\[15pt]
\displaystyle\mbox{}\hskip-0.4cm
\sum_{\alpha,\beta=1}^3\int_{R^3}|D^{\sigma}P_ju|^{k-2}D^{\sigma}P_ju_\beta\partial_{\alpha}D^{\sigma}P_j(u_\alpha
u_\beta)=
\\[15pt]
\displaystyle
-(k-2)\sum_{\alpha,\beta=1}^3\int_{R^3}D^{\sigma}P_j(u_\alpha
u_\beta)|D^{\sigma}P_ju|^{k-3}\partial_{\alpha}|D^{\sigma}P_ju|D^{\sigma}P_ju_\beta\\[15pt]
\displaystyle -\sum_{\alpha,\beta=1}^3
\int_{R^3}D^{\sigma}P_j(u_\alpha
u_\beta)|D^{\sigma}P_ju|^{k-2}D^{\sigma}P_j\partial_{\alpha}u_\beta\end{array}\eqno{(5.9)}$$
and
$$\begin{array}{l}\displaystyle\mbox{}\hskip-1cm\int_{R^3}|D^{\sigma}P_ju|^{k-2}D^{\sigma}P_juD^{\sigma}P_j\nabla p
\\[15pt]
\displaystyle
=-\sum_{\alpha}^3\int_{R^3}\partial_{\alpha}\left(|D^{\sigma}P_ju|^{k-2}D^{\sigma}P_ju_{\alpha}\right)D^{\sigma}P_jp\\[15pt]
\displaystyle=
-(k-2)\sum_{\alpha}^3\int_{R^3}|D^{\sigma}P_ju|^{k-3}D^{\sigma}P_ju_\alpha\partial_{\alpha}|D^{\sigma}P_ju|D^{\sigma}P_jp
,\end{array}\eqno{(5.10)}$$
where $u=(u_1,u_2,u_3)$ and $\mbox{div}u=0$ is used.

Plug (5.7)-(5.10) into (5.6) and we have
$$\begin{array}{l}\displaystyle\frac1{k}\frac{d}{dt}\int_{R^3}|D^{\sigma}P_ju|^{k}+
\nu\int_{R^3}|D^{\sigma}P_ju|^{k-2}|D^{\sigma}P_j\nabla u|^2
\\[15pt]\displaystyle\mbox{}\hskip0.5cm
+
(k-2)\nu\int_{R^3}|D^{\sigma}P_ju|^{k-2}|\nabla|D^{\sigma}P_ju||^2
\\[15pt]
\displaystyle=
(k-2)\sum_{\alpha,\beta=1}^3\int_{R^3}D^{\sigma}P_j(u_\alpha
u_\beta)|D^{\sigma}P_ju|^{k-3}\partial_{\alpha}|D^{\sigma}P_ju|D^{\sigma}P_ju_\beta\\[15pt]
\displaystyle\mbox{}\hskip0.5cm +\sum_{\alpha,\beta=1}^3
\int_{R^3}D^{\sigma}P_j(u_\alpha
u_\beta)|D^{\sigma}P_ju|^{k-2}D^{\sigma}P_j\partial_{\alpha}u_\beta\\[15pt]
\displaystyle\mbox{}\hskip0.5cm
+(k-2)\sum_{\alpha,\beta=1}^3\int_{R^3}|D^{\sigma}P_ju|^{k-3}D^{\sigma}P_ju_\alpha\partial_{\alpha}|D^{\sigma}P_ju|D^{\sigma}P_jp
\\[15pt]
\displaystyle:=(k-2)I_1+I_2+(k-2)I_3.\end{array}\eqno{(5.11)}$$

\bigskip

{\it Step 2. Estimates of $I_1,I_2$ and $I_3$ in (5.11).}\hsa
From H\"{o}lder's inequality,
$$\begin{array}{l}\displaystyle|I_1|\leq\int_{R^3}|D^{\sigma}P_j(u\otimes u)||D^{\sigma}P_ju|^{k-2}|\nabla|D^{\sigma}P_ju||\\[15pt]
\displaystyle\mbox{}\hskip0.5cm\leq
||D^{\sigma}P_j(u\otimes u)||_{\frac{3k}{k+1}}
\left(||D^{\sigma}P_ju||_{k}^{k}\right)^{\frac{k-2}{12k}}
\left(||D^{\sigma}P_ju||_{5k}^{5k}\right)^{\frac{k-2}{12k}}
\\[15pt]\displaystyle\mbox{}\hskip1cm\times
\left(\int_{R^3}|D^{\sigma}P_ju|^{k-2}|\nabla|D^{\sigma}P_ju||^2\right)^{\frac12}.
\end{array}$$
From Lemma 4.5, it follows that
$$\begin{array}{l}\displaystyle\mbox{}\hskip-0cm
||D^{\sigma}P_ju||_{k}^{k}
\leq Ck^22^{-2j}
\int_{R^3}|D^{\sigma}P_ju|^{k-2}|\nabla D^{\sigma}P_ju|^2.
\end{array}$$
And then
$$\begin{array}{l}\displaystyle|I_1|\leq
Ck^{\frac{k-2}{6k}}2^{-\frac{k-2}{6k}j}
||D^{\sigma}P_j(u\otimes u)||_{\frac{3k}{k+1}}
\left(||D^{\sigma}P_ju||_{5k}^{5k}\right)^{\frac{k-2}{12k}}
\\[15pt]\displaystyle\mbox{}\hskip1cm\times
\left(
\int_{R^3}|D^{\sigma}P_ju|^{k-2}|\nabla D^{\sigma}P_ju|^2\right)^{\frac{k-2}{12k}}
\left(\int_{R^3}|D^{\sigma}P_ju|^{k-2}|\nabla|D^{\sigma}P_ju||^2\right)^{\frac12}.
\end{array}$$

By Young's inequality,
$$\begin{array}{l}\displaystyle|I_1|\leq
Ck^{\frac{3(k-2)}{5k+2}}2^{-\frac{2(k-2)}{5k+2}j}\nu^{-\frac{7k-2}{5k+2}}
||D^{\sigma}P_j(u\otimes u)||_{\frac{3k}{k+1}}^{\frac{12k}{5k+2}}
\left(||D^{\sigma}P_ju||_{5k}^{5k}\right)^{\frac{k-2}{5k+2}}
\\[15pt]\displaystyle\mbox{}\hskip0.6cm
+\frac{\nu}{3(k-2)}
\int_{R^3}|D^{\sigma}P_ju|^{k-2}|\nabla D^{\sigma}P_ju|^2
+\frac{\nu}{3}
\int_{R^3}|D^{\sigma}P_ju|^{k-2}|\nabla|D^{\sigma}P_ju||^2
\\[15pt]\displaystyle\mbox{}\hskip0.2cm
\leq
Ck^{2}2^{-\frac{j}4}\left(1+\nu^{-2}\right)
||D^{\sigma}P_j(u\otimes u)||_{\frac{3k}{k+1}}^{\frac{12k}{5k+2}}
\left(||D^{\sigma}P_ju||_{5k}^{5k}\right)^{\frac{k-2}{5k+2}}
\\[15pt]\displaystyle\mbox{}\hskip0.6cm
+\frac{\nu}{3(k-2)}
\int_{R^3}|D^{\sigma}P_ju|^{k-2}|\nabla D^{\sigma}P_ju|^2
+\frac{\nu}{3}
\int_{R^3}|D^{\sigma}P_ju|^{k-2}|\nabla|D^{\sigma}P_ju||^2.
\end{array}\eqno{(5.12)}$$

We estimate $I_2$ as almost same as $I_1$.
Actually, we only need replace the term
$\int_{R^3}|D^{\sigma}P_ju|^{k-2}|\nabla|D^{\sigma}P_ju||^2$ by the term
$\int_{R^3}|D^{\sigma}P_ju|^{k-2}|D^{\sigma}P_j\nabla u|^2$ in the estimate of $I_1$.
That is, we have
$$\begin{array}{l}\displaystyle|I_2|
\leq
Ck^{2}2^{-\frac{j}4}\left(1+\nu^{-2}\right)
||D^{\sigma}P_j(u\otimes u)||_{\frac{3k}{k+1}}^{\frac{12k}{5k+2}}
\left(||D^{\sigma}P_ju||_{5k}^{5k}\right)^{\frac{k-2}{5k+2}}
\\[15pt]\displaystyle\mbox{}\hskip1cm
+\frac{\nu}{3}
\int_{R^3}|D^{\sigma}P_ju|^{k-2}|\nabla D^{\sigma}P_ju|^2.
\end{array}\eqno{(5.13)}$$

$I_3$ can be estimated also by the similar way. From H\"{o}lder's inequality,
$$\begin{array}{l}\displaystyle|I_3|\leq
||D^{\sigma}P_jp||_{\frac{3k}{k+1}}
\left(||D^{\sigma}P_ju||_{k}^{k}\right)^{\frac{k-2}{12k}}
\left(||D^{\sigma}P_ju||_{5k}^{5k}\right)^{\frac{k-2}{12k}}
\\[15pt]\displaystyle\mbox{}\hskip1cm\times
\left(\int_{R^3}|D^{\sigma}P_ju|^{k-2}|\nabla|D^{\sigma}P_ju||^2\right)^{\frac12}.
\end{array}$$
From Lemma 3.7,
$$||D^{\sigma}P_jp||_{\frac{3k}{k+1}}\leq C||D^{\sigma}P_j(u\otimes u)||_{\frac{3k}{k+1}},$$
where $C$ is a universal constant.
Therefore
$$\begin{array}{l}\displaystyle|I_3|
\leq
Ck^{2}2^{-\frac{j}4}\left(1+\nu^{-2}\right)
||D^{\sigma}P_j(u\otimes u)||_{\frac{3k}{k+1}}^{\frac{12k}{5k+2}}
\left(||D^{\sigma}P_ju||_{5k}^{5k}\right)^{\frac{k-2}{5k+2}}
\\[15pt]\displaystyle\mbox{}\hskip0.6cm
+\frac{\nu}{3(k-2)}
\int_{R^3}|D^{\sigma}P_ju|^{k-2}|\nabla D^{\sigma}P_ju|^2
+\frac{\nu}{3}
\int_{R^3}|D^{\sigma}P_ju|^{k-2}|\nabla|D^{\sigma}P_ju||^2.
\end{array}\eqno{(5.14)}$$

Plug (5.12)-(5.14) into (5.11). We arrive at
$$\begin{array}{l}\displaystyle\mbox{}\hskip-0cm\frac{d}{dt}||D^{\sigma}P_j
u||_{k}^{k}
\leq
Ck^{4}2^{-\frac{j}4}\left(1+\nu^{-2}\right)
\\[15pt]\displaystyle\mbox{}\hskip2cm\times
||D^{\sigma}P_j(u\otimes u)||_{\frac{3k}{k+1}}^{\frac{12k}{5k+2}}
\left(||D^{\sigma}P_ju||_{5k}^{5k}\right)^{\frac{k-2}{5k+2}}.
\end{array}\eqno{(5.15)}$$

\bigskip

{\it Step 3. Estimate of $||D^{\sigma}P_j(u\otimes u)||_{\frac{3k}{k+1}}^{\frac{12k}{5k+2}}$ in (5.15).}\hsa
In Corollary 2.8, we set $\hat q=\frac{12k}{5k+2}$, $q=\frac{3k}{k+1}$, $q_0=\frac{3k}{k-2}$ and $q_1=k$ and then using
it, we have
$$\begin{array}{l}\displaystyle\mbox{}\hskip0cm
||D^{\sigma}P_j(u\otimes u)||_{\frac{3k}{k+1}}^{\frac{12k}{5k+2}}
\leq
C\alpha_j||u||_2^{\frac{24k}{5k+2}}
\\[15pt]\displaystyle\mbox{}\hskip2cm
+
C\sum_{m=\max\{j_0,j-2\}}^{\infty}2^{\sigma(j-m)}||D^{\sigma}P_mu||_{k}^{\frac{12k}{5k+2}}||u||_{\frac{3k}{k-2}}^{\frac{12k}{5k+2}}
\\[15pt]\displaystyle\mbox{}\hskip0cm
\leq C\left(||u||_2^{\frac{24k}{5k+2}}+||u||_{\frac{3k}{k-2}}^{\frac{12k}{5k+2}}\right)
\left(\alpha_j+\sum_{m=\max\{j_0,j-2\}}^{\infty}2^{\sigma(j-m)}||D^{\sigma}P_mu||_{k}^{\frac{12k}{5k+2}}\right)
,
\end{array}$$
where $\alpha_j$ is given by (2.3) and $C$ are universal constants.
Since $$||u||_{\frac{3k}{k-2}}\leq||u||_2^{\lambda}||u||_4^{1-\lambda}\leq||u||_2+||u||_4,$$
where $\frac{k-2}{3k}=\frac{\lambda}2+\frac{1-\lambda}4$, it is easy to see that
$$||u||_2^{\frac{24k}{5k+2}}+||u||_{\frac{3k}{k-2}}^{\frac{12k}{5k+2}}\leq C
\left(1+||u||_2^5+||u||_4^{\frac83}\right).$$
For simplicity, we denote
$$\cc_0(t):=1+||u||_2^5+||u||_4^{\frac83}.\eqno{(5.16)}$$
Then
$$\begin{array}{l}\displaystyle\mbox{}\hskip0cm
||D^{\sigma}P_j(u\otimes u)||_{\frac{3k}{k+1}}^{\frac{12k}{5k+2}}
\leq C\cc_0(t)
\left(\alpha_j+\sum_{m=\max\{j_0,j-2\}}^{\infty}2^{\sigma(j-m)}||D^{\sigma}P_mu||_{k}^{\frac{12k}{5k+2}}\right)
.
\end{array}$$
From (5.15), it follows that
$$\begin{array}{l}\displaystyle\mbox{}\hskip-0cm\frac{d}{dt}||D^{\sigma}P_ju||_{k}^{k}
\leq
C\left(1+\nu^{-2}\right)\cc_0(t)k^{4}2^{-\frac{j}4}\Bigg\{
\alpha_j\left(||D^{\sigma}P_ju||_{5k}^{5k}\right)^{\frac{k-2}{5k+2}}
\\[15pt]\displaystyle\mbox{}\hskip0.8cm
+
\sum_{m=\max\{j_0,j-2\}}^{\infty}2^{\sigma(j-m)}||D^{\sigma}P_mu||_{k}^{\frac{12k}{5k+2}}
\left(||D^{\sigma}P_ju||_{5k}^{5k}\right)^{\frac{k-2}{5k+2}}
\Bigg\}.
\end{array}$$
From Young's inequality, we see
$$\begin{array}{l}\displaystyle
||D^{\sigma}P_mu||_{k}^{\frac{12k}{5k+2}}
\left(||D^{\sigma}P_ju||_{5k}^{5k}\right)^{\frac{k-2}{5k+2}}
=
\left(||D^{\sigma}P_mu||_{k}^k\right)^{\frac{12}{5k+2}}
\left(\left(||D^{\sigma}P_ju||_{5k}^{5k}\right)^{\frac15}\right)^{\frac{5k-10}{5k+2}}
\\[15pt]\displaystyle\mbox{}\hskip0.8cm
\leq C\left(\frac1{k^4}||D^{\sigma}P_mu||_{k}^k\right)^{\frac{12}{5k+2}}
\left(\left(||D^{\sigma}P_ju||_{5k}^{5k}\right)^{\frac15}\right)^{\frac{5k-10}{5k+2}}
\\[15pt]\displaystyle\mbox{}\hskip0.8cm
\leq C\left(\frac1{k^4}||D^{\sigma}P_mu||_{k}^k+
\left(||D^{\sigma}P_ju||_{5k}^{5k}\right)^{\frac15}\right).
\end{array}$$
It follows that
$$\begin{array}{l}\displaystyle\mbox{}\hskip-0cm
\sum_{m=\max\{j_0,j-2\}}^{\infty}2^{\sigma(j-m)}||D^{\sigma}P_mu||_{k}^{\frac{12k}{5k+2}}
\left(||D^{\sigma}P_ju||_{5k}^{5k}\right)^{\frac{k-2}{5k+2}}
\\[15pt]\displaystyle\mbox{}\hskip0.8cm
\leq C\sum_{m=\max\{j_0,j-2\}}^{\infty}2^{\sigma(j-m)}
\left(\frac1{k^4}||D^{\sigma}P_mu||_{k}^k+
\left(||D^{\sigma}P_ju||_{5k}^{5k}\right)^{\frac15}\right)
\\[15pt]\displaystyle\mbox{}\hskip0.8cm
\leq C\left(\left(||D^{\sigma}P_ju||_{5k}^{5k}\right)^{\frac15}+\sum_{m=\max\{j_0,j-2\}}^{\infty}2^{\sigma(j-m)}
\frac1{k^4}||D^{\sigma}P_mu||_{k}^k\right)
.
\end{array}$$
It is clear that
$$\alpha_j\left(||D^{\sigma}P_ju||_{5k}^{5k}\right)^{\frac{k-2}{5k+2}}\leq\alpha_j
\left(1+\left(||D^{\sigma}P_ju||_{5k}^{5k}\right)^{\frac15}\right).$$
Therefore
$$\begin{array}{l}\displaystyle\mbox{}\hskip-0cm\frac{d}{dt}||D^{\sigma}P_ju||_{k}^{k}
\leq
C\left(1+\nu^{-2}\right)\cc_0(t)k^{4}2^{-\frac{j}4}
\\[15pt]\displaystyle\mbox{}\hskip0.5cm\times
\Bigg\{
\alpha_j+\left(||D^{\sigma}P_ju||_{5k}^{5k}\right)^{\frac15}
+
\sum_{m=\max\{j_0,j-2\}}^{\infty}\frac{2^{\sigma(j-m)}}{k^4}||D^{\sigma}P_mu||_{k}^{k}
\Bigg\}
.
\end{array}\eqno{(5.17)}$$

\bigskip

{\it Step 4. Taking sum of $j$.}\hsa
Take sum of $j$ from $j_0=1$ to $\infty$ on both sides of (5.17). Since
$$\begin{array}{l}\displaystyle
\sum_{j=j_0}^{\infty}\sum_{m=\max\{j_0,j-2\}}^{\infty}\frac{2^{\sigma(j-m)}}{k^4}||D^{\sigma}P_mu||_{k}^{k}
=
\sum_{j=3}^{\infty}\sum_{m=j-2}^{\infty}\frac{2^{\sigma(j-m)}}{k^4}||D^{\sigma}P_mu||_{k}^{k}
\\[15pt]\displaystyle\mbox{}\hskip2cm
+
\sum_{m=j_0}^{\infty}\frac{2^{\sigma(1-m)}}{k^4}||D^{\sigma}P_mu||_{k}^{k}
+
\sum_{m=j_0}^{\infty}\frac{2^{\sigma(2-m)}}{k^4}||D^{\sigma}P_mu||_{k}^{k}
\\[15pt]\displaystyle\mbox{}\hskip1cm
\leq C\sum_{m=j_0}^{\infty}\frac{1}{k^4}||D^{\sigma}P_mu||_{k}^{k},
\end{array}$$
by Lemma 4.2 and 4.6, we have
$$\begin{array}{l}\displaystyle\mbox{}\hskip-0cm\frac{d}{dt}\sum_{j=j_0}^{\infty}||D^{\sigma}P_ju||_{k}^{k}
\leq
C\left(1+\nu^{-2}\right)\cc_0(t)k^{4}
\\[15pt]\displaystyle\mbox{}\hskip1cm
\times
\Bigg\{
1+\sum_{j=j_0}^{\infty}2^{-\frac{j}4}\left(||D^{\sigma}P_ju||_{5k}^{5k}\right)^{\frac15}

+\sum_{j=j_0}^{\infty}\frac{1}{k^4}||D^{\sigma}P_ju||_{k}^{k}
\Bigg\}
.
\end{array}\eqno{(5.18)}$$

\bigskip

{\it Step 5. Dividing (5.18) by $2^{B_k}$ and taking sum of $k$.}\hsa
Divide by $2^{B_k}$ on both sides of (5.18) and consequently,
$$\begin{array}{l}\displaystyle\mbox{}\hskip-0cm
\frac{d}{dt}\sum_{j=1}^{\infty}
\frac{||D^{\sigma}P_ju||_{k}^{k}}{2^{B_k}}
\leq
C\left(1+\nu^{-2}\right)\cc_0(t)k^{4}
\\[15pt]\displaystyle\mbox{}\hskip0cm\times
\Bigg\{\frac{1}{2^{B_k}}+
\sum_{j=1}^{\infty}\frac{1}{k^4}\frac{||D^{\sigma}P_ju||_{k}^{k}}{2^{B_{k}}}
+\frac{2^{B_{5k}/5}}{2^{B_k}}
\sum_{j=1}^{\infty}2^{-\frac{j}4}
\left(\frac{||D^{\sigma}P_ju||_{5k}^{5k}}{2^{B_{5k}}}\right)^{\frac15}\Bigg\}.
\end{array}\eqno{(5.19)}$$
From (5.2), we see that
$$\begin{array}{l}\displaystyle\mbox{}\hskip-0cm
\frac{2^{B_{5k}/5}}{2^{B_k}}=
2^{\left(\frac{1}{\sqrt{5}}-1\right)\sqrt{k}}
\leq\frac{C}{k^6}
\end{array}$$
for some universal constant $C$.
It is clear that
$$\begin{array}{l}\displaystyle\mbox{}\hskip-0cm
\frac{1}{2^{B_k}}\leq\frac{C}{k^6}
\end{array}$$
for some universal constant $C$.
From (5.19), it follows that
$$\begin{array}{l}\displaystyle\mbox{}\hskip-0cm
\frac{d}{dt}\sum_{j=1}^{\infty}
\frac{||D^{\sigma}P_ju||_{k}^{k}}{2^{B_k}}
\leq
C\left(1+\nu^{-2}\right)\cc_0(t)
\\[15pt]\displaystyle\mbox{}\hskip0cm\times
\Bigg\{\frac{1}{k^2}+
\sum_{j=1}^{\infty}\frac{||D^{\sigma}P_ju||_{k}^{k}}{2^{B_{k}}}
+
\sum_{j=1}^{\infty}\frac{1}{2^{\frac{j}4}k^2}
\left(\frac{||D^{\sigma}P_ju||_{5k}^{5k}}{2^{B_{5k}}}\right)^{\frac15}\Bigg\}.
\end{array}\eqno{(5.20)}$$

Take sum of $k$ from $k_0$ to $\infty$ on both sides of (5.20). From (5.3) and Lemma 4.2, we obtain
$$\begin{array}{l}\displaystyle\mbox{}\hskip-0cm
\frac{d}{dt}\sum_{k=k_0}^{\infty}\sum_{j=1}^{\infty}
\frac{||D^{\sigma}P_ju||_{k}^{k}}{2^{B_k}}
\leq
C\left(1+\nu^{-2}\right)\cc_0(t)
\\[15pt]\displaystyle\mbox{}\hskip0cm\times
\Bigg\{\sum_{k=k_0}^{\infty}\frac{1}{k^2}+
\sum_{k=k_0}^{\infty}\sum_{j=1}^{\infty}\frac{||D^{\sigma}P_ju||_{k}^{k}}{2^{B_{k}}}
+
\sum_{k=k_0}^{\infty}\sum_{j=1}^{\infty}\frac{1}{2^{\frac{j}4}k^2}
\left(\frac{||D^{\sigma}P_ju||_{5k}^{5k}}{2^{B_{5k}}}\right)^{\frac15}\Bigg\}.
\end{array}$$

It is clear that
$$\sum_{k=k_0}^{\infty}\frac{1}{k^2}\leq C$$
and from Young's inequality,
$$\begin{array}{l}\displaystyle\mbox{}\hskip-0cm
\sum_{k=k_0}^{\infty}\sum_{j=1}^{\infty}\frac{1}{2^{\frac{j}4}k^2}
\left(\frac{||D^{\sigma}P_ju||_{5k}^{5k}}{2^{B_{5k}}}\right)^{\frac15}
\leq C\Bigg\{1+\sum_{k=k_0}^{\infty}\sum_{j=1}^{\infty}
\frac{||D^{\sigma}P_ju||_{5k}^{5k}}{2^{B_{5k}}}\Bigg\}.
\end{array}$$
Therefore we have
$$\begin{array}{l}\displaystyle\mbox{}\hskip-0cm
\frac{d}{dt}\sum_{k=k_0}^{\infty}\sum_{j=1}^{\infty}
\frac{||D^{\sigma}P_ju||_{k}^{k}}{2^{B_k}}
\leq
C\left(1+\nu^{-2}\right)\cc_0(t)
\Bigg\{1+\sum_{k=k_0}^{\infty}\sum_{j=1}^{\infty}
\frac{||D^{\sigma}P_ju||_{k}^{k}}{2^{B_k}}\Bigg\}.
\end{array}\eqno{(5.21)}$$

\bigskip

{\it Step 6. Proof of (5.4).}\hsa
From (5.21), it follows that
$$\begin{array}{l}\displaystyle\mbox{}\hskip-0cm
\frac{d}{dt}\Bigg\{1+\sum_{k=k_0}^{\infty}\sum_{j=1}^{\infty}
\frac{||D^{\sigma}P_ju||_{k}^{k}}{2^{B_k}}\Bigg\}
\leq
C\left(1+\nu^{-2}\right)\cc_0(t)
\Bigg\{1+\sum_{k=k_0}^{\infty}\sum_{j=1}^{\infty}
\frac{||D^{\sigma}P_ju||_{k}^{k}}{2^{B_k}}\Bigg\}.
\end{array}$$
From Gronwall's inequality, we have
$$\begin{array}{l}\displaystyle
\sum_{k=k_0}^{\infty}\sum_{j=j_0}^{\infty}
\frac{||\dd{\sigma}P_j
u(t)||_{k}^{k}}{2^{B_k}}
\leq
e^{C\left(1+\nu^{-2}\right)\int_0^T\cc_0(t))ds}-1
\\[15pt]\displaystyle\mbox{}\hskip1cm+
e^{C\left(1+\nu^{-2}\right)\int_0^T\cc_0(t)ds}
\sum_{k=k_0}^{\infty}\sum_{j=j_0}^{\infty}\frac{||\dd{\sigma}P_j
u_0||_{k}^{k}}{2^{B_k}}
\end{array}$$
for any $0\leq t\leq T'$.
In view of (5.16) and Corollary 3.2, we have
$$
\begin{array}{l}\displaystyle
\int_0^T\cc_0(t)dt=
\int_0^T\left(1+||u||_2^{5}+||u||_{4}^{\frac83}\right)dt
\leq
\left(1+||u_0||_2^5\right)T+\frac{C}{\nu}||u_0||_2^{\frac83}.
\end{array}
$$
Then (5.4) follows clearly.
$\Box$

\bigskip

\section{Low frequency part}

From Theorem 5.1 ({\it uniform bound estimate}) to derive our new a priori estimate, we only need to remove Condition (5.3).
To do this, we separate the series (5.1) into two parts, low frequency part (finite $j$) and high frequency part,
and show the convergence of them respectively.
In this section, we will study the low frequency part which is much more difficult than the high frequency part
and the result, Theorem 6.1 is the second key step to approach our new a priori estimates.
We design the following series
$$\displaystyle\sum_{k=k_0}^{\infty}\sum_{j=j_0}^{J_0}\fd{||D^{\sigma}P_{j}
u(t)||_{{k}}^{{k}}}{2^{\hat B_{k}}},\eqno{(6.1)}$$
where $j_0=1$, $k_0=100$, $J_0$ will be given later and
the power $$\hat B_k=\left(B-\frac{1}{\sqrt{k}}\right)k+2^B.\eqno{(6.2)}$$
Note here the constant $B$ will be chosen the same as in (5.2).
We will see that (6.1) can not blow up before any given time $T$ if $B$ is large enough.

Our main theorem of this section is:

\bigskip

\theorem{6.1}{Let $j_0=1$, $k_0=100$, $\sigma=2$, $T>0$ and $B>1$ be real numbers,
$\hat B_k$ be given by (6.2) for any $k\geq k_0$, and $u_0$ be a function satisfying (1.2).
Suppose $u$ and $p$ satisfy Condition (S).
Let $$J_0=\left[\frac{8B}{\sigma}\right].\eqno{(6.3)}$$
There exists $\tilde B_1>0$ depending only on $\nu$, $T$ and $u_0$ such that
if $B\geq\tilde B_1$, then
$$\displaystyle\sum_{k=k_0}^{\infty}\sum_{j=j_0}^{J_0}\fd{||D^{\sigma}P_{j}
u(t)||_{{k}}^{{k}}}{2^{\hat B_{k}}}\leq1\eqno{(6.4)}$$
for any $t\in[0,T]$.
}

\bigskip

We establish Theorem 6.1 by the following two lemmas, where $B>1$ will be determined later.

\bigskip

\lemma{6.2}{Suppose all the assumptions of Theorem 6.1 hold. Then
$$\begin{array}{l}\displaystyle\mbox{}\hskip-0cm\frac{d}{dt}\sum_{k=k_0}^{2k_0-1}\sum_{j=j_0}^{J_0}
\frac{||D^{\sigma}P_ju||_{k}^{k}}{2^{\hat B_k}}
\leq
\frac{ C}{\nu}||u||_2^4\left(2^{-\frac14B}+
\sum_{k=k_0}^{2k_0-1}\sum_{j=j_0}^{J_0}\frac{||D^{\sigma}P_ju||_{k}^{k}}{2^{\hat B_{k}}}\right).
\end{array}\eqno{(6.5)}$$
for $0\leq t\leq T$.
}

\proof Let $k_0\leq k<2k_0$, $j_0\leq j\leq J_0$ and $0\leq t\leq T$. We divide the proof into four steps.

\bigskip

{\it  Step 1.}\hsa
By the same arguments to derive (5.11),
we have
$$\begin{array}{l}\displaystyle\frac1{k}\frac{d}{dt}\int_{R^3}|D^{\sigma} P_ju|^{k}+
\nu\int_{R^3}|D^{\sigma} P_ju|^{k-2}|\nabla D^{\sigma}P_ju|^2
\\[15pt]\displaystyle\mbox{}\hskip1.5cm
+
(k-2)\nu\int_{R^3}|D^{\sigma} P_ju|^{k-2}|\nabla|D^{\sigma} P_j u||^2
\\[15pt]\displaystyle\mbox{}\hskip0.2cm
\leq
(k-2)\int_{R^3}|D^{\sigma} P_ju|^{k-2}|\nabla |D^{\sigma} P_j u|| |D^{\sigma} P_j\left(u\otimes u\right)|
\\[15pt]\displaystyle\mbox{}\hskip2cm
+
\int_{R^3}|D^{\sigma} P_ju|^{k-2}|\nabla D^{\sigma} P_j u| |D^{\sigma} P_j\left(u\otimes u\right)|
\\[15pt]\displaystyle\mbox{}\hskip1.5cm
+
(k-2)\int_{R^3}|D^{\sigma} P_ju|^{k-2}|\nabla |D^{\sigma} P_j u|| |D^{\sigma}P_jp|
\\[15pt]\displaystyle\mbox{}\hskip0.2cm
:=(k-2)I_1+I_2+(k-2)I_3
.\end{array}\eqno{(6.6)}$$

\bigskip

{\it Step 2. Estimates of $I_1,I_2$ and $I_3$ in (6.6).}\hsa
From H\"{o}lder's inequality,
$$\begin{array}{l}\displaystyle |I_1|\leq
||D^{\sigma} P_j\left(u\otimes u\right)||_k
||D^{\sigma} P_ju||_k^{\frac{k-2}2}
\left(\int_{R^3}|D^{\sigma} P_j u|^{k-2}|\nabla |D^{\sigma}P_j u||^2\right)^{\frac12}
\\[15pt]\displaystyle\mbox{}\hskip0.5cm
\leq
\frac4{\nu}
||D^{\sigma} P_j\left(u\otimes u\right)||_k^2
||D^{\sigma} P_ju||_k^{k-2}
\\[15pt]\displaystyle\mbox{}\hskip2cm
+
\frac{\nu}4
\int_{R^3}|D^{\sigma} P_j u|^{k-2}|\nabla|D^{\sigma}P_j u||^2.
\end{array}$$
and
$$\begin{array}{l}\displaystyle |I_2|\leq
||D^{\sigma} P_j\left(u\otimes u\right)||_k
||D^{\sigma} P_ju||_k^{\frac{k-2}2}
\left(\int_{R^3}|D^{\sigma} P_j u|^{k-2}|\nabla D^{\sigma}P_j u|^2\right)^{\frac12}
\\[15pt]\displaystyle\mbox{}\hskip0.5cm
\leq
\frac1{\nu}
||D^{\sigma} P_j\left(u\otimes u\right)||_k^2
||D^{\sigma} P_ju||_k^{k-2}
\\[15pt]\displaystyle\mbox{}\hskip2cm
+
\nu
\int_{R^3}|D^{\sigma} P_j u|^{k-2}|\nabla D^{\sigma}P_j u|^2.
\end{array}$$

From Lemma 3.7, we have
$$||D^{\sigma} P_jp||_k\leq Ck||D^{\sigma} P_j\left(u\otimes u\right)||_k,$$
where $C$ is a universal constant. Using this to estimate $I_3$, we obtain
$$\begin{array}{l}\displaystyle |I_3|\leq
||D^{\sigma} P_jp||_k
||D^{\sigma} P_ju||_k^{\frac{k-2}2}
\left(\int_{R^3}|D^{\sigma} P_j u|^{k-2}|\nabla |D^{\sigma}P_j u||^2\right)^{\frac12}
\\[15pt]\displaystyle\mbox{}\hskip0.5cm
\leq
Ck||D^{\sigma} P_j\left(u\otimes u\right)||_k
||D^{\sigma} P_ju||_k^{\frac{k-2}2}
\left(\int_{R^3}|D^{\sigma} P_j u|^{k-2}|\nabla |D^{\sigma}P_j u||^2\right)^{\frac12}
\\[15pt]\displaystyle\mbox{}\hskip0.5cm
\leq
\frac{Ck^2}{\nu}
||D^{\sigma} P_j\left(u\otimes u\right)||_k^2
||D^{\sigma} P_ju||_k^{k-2}
\\[15pt]\displaystyle\mbox{}\hskip2.8cm
+
\frac{\nu}4
\int_{R^3}|D^{\sigma} P_j u|^{k-2}|\nabla|D^{\sigma}P_j u||^2.
\end{array}$$

Plug the estimates of $I_1,I_2$ and $I_3$  into (6.6) and we arrive at
$$\begin{array}{l}\displaystyle\mbox{}\hskip0cm
\frac{d}{dt}||D^{\sigma} P_j u||_{k}^{k}+
\frac{\nu}2k(k-2)
\int_{R^3}|D^{\sigma} P_j u|^{k-2}|\nabla|D^{\sigma}P_j u||^2
\\[15pt]\displaystyle\mbox{}\hskip1cm
\leq
\frac{ Ck^4}{\nu}||D^{\sigma} P_j\left(u\otimes u\right)||_k^2
||D^{\sigma} P_ju||_k^{k-2}.
\end{array}\eqno{(6.7)}$$

\bigskip

{\it Step 3. Simplifying (6.7).}\hsa
By Lemma 2.1, 2.3 and 2.6, we have
$$\begin{array}{l}\displaystyle\mbox{}\hskip0cm
||D^{\sigma} P_j\left(u\otimes u\right)||_k
\leq C2^{(\sigma+3-\frac3k)j}||P_j\left(u\otimes u\right)||_1
\leq C2^{(\sigma+3)J_0}||u||_2^2,
\end{array}$$
where $j\leq J_0$ is used. From $k_0\leq k<2k_0$, (6.2) and (6.3), we deduce
$$
2^{(\sigma+3)J_0}\leq2^{\frac{8(\sigma+3)}{\sigma}B}\leq C2^{(2^B-B)/k}\leq C2^{(\hat B_k-B)/k}
$$
for any $B>1$, where $C$ is a universal constant. It follows that
$$\begin{array}{l}\displaystyle\mbox{}\hskip0cm
||D^{\sigma} P_j\left(u\otimes u\right)||_k
\leq C2^{(\hat B_k-B)/k}||u||_2^2,
\end{array}$$
as $j\leq J_0$ and $k_0\leq k<2k_0$.

Plug this estimate into (6.7) and we have
$$\begin{array}{l}\displaystyle\mbox{}\hskip0cm
\frac{d}{dt}||D^{\sigma} P_j u||_{k}^{k}
\leq
\frac{ C}{\nu}2^{2(\hat B_k-B)/k}||u||_2^4
||D^{\sigma} P_ju||_k^{k-2},
\end{array}\eqno{(6.8)}$$
where the term
$\frac{\nu}2k(k-2)
\int_{R^3}|D^{\sigma} P_j u|^{k-2}|\nabla|D^{\sigma}P_j u||^2$ is omitted and $k_0\leq k<2k_0$ is used.

\bigskip

{\it Step 4. Proof of (6.5).}\hsa
Divide $2^{\hat B_k}$ on the both sides of (6.8) and then
$$\begin{array}{l}\displaystyle\mbox{}\hskip-0cm\frac{d}{dt}
\frac{||D^{\sigma}P_ju||_{k}^{k}}{2^{\hat B_k}}
\leq
\frac{ C}{\nu}||u||_2^4\left(2^{-2B/k}
\left(\frac{||D^{\sigma}P_ju||_{k}^{k}}{2^{\hat B_{k}}}\right)^{\frac{k-2}k}\right).
\end{array}$$
From Young's inequality, we have
$$\begin{array}{l}\displaystyle\mbox{}\hskip-0cm\frac{d}{dt}
\frac{||D^{\sigma}P_ju||_{k}^{k}}{2^{\hat B_k}}
\leq
\frac{ C}{\nu}||u||_2^4\left(2^{-B}+
\frac{||D^{\sigma}P_ju||_{k}^{k}}{2^{\hat B_{k}}}\right).
\end{array}\eqno{(6.9)}$$

Take sum of $j$ from $j_0$ to $J_0$ and of $k$ from $k_0$ to $2k_0-1$ on both sides of (6.9).
We have
$$\begin{array}{l}\displaystyle\mbox{}\hskip-0cm\frac{d}{dt}\sum_{k=k_0}^{2k_0-1}\sum_{j=j_0}^{J_0}
\frac{||D^{\sigma}P_ju||_{k}^{k}}{2^{\hat B_k}}
\leq
\frac{ C}{\nu}||u||_2^4\left(\sum_{k=k_0}^{2k_0-1}\sum_{j=j_0}^{J_0}2^{-B}+
\sum_{k=k_0}^{2k_0-1}\sum_{j=j_0}^{J_0}\frac{||D^{\sigma}P_ju||_{k}^{k}}{2^{\hat B_{k}}}\right).
\end{array}$$
From (6.3), it follows that
$$\begin{array}{l}\displaystyle
\sum_{k=k_0}^{2k_0-1}\sum_{j=j_0}^{J_0}2^{-B}=k_0J_02^{-B}\leq k_0\frac{8}{\sigma}B2^{-B}\leq C2^{-\frac14B}
\end{array}$$
for any $B>1$ where $C$ is a universal constant. Then we see (6.5) clearly.
$\Box$

\bigskip

\lemma{6.3}{Suppose all the assumptions of Theorem 6.1 hold. Then for any $k'\geq2k_0$, we have
$$\begin{array}{l}\displaystyle\mbox{}\hskip0cm
\frac{d}{dt}\sum_{k=2k_0}^{k'}\sum_{j=j_0}^{J_0}\frac{||D^{\sigma} P_j u||_{k}^{k}}{2^{\hat B_k}}
\leq
C\left(1+\frac1{\nu}\right)^{4}\left(1+||u||_2^2\right)^5
\\[15pt]\displaystyle\mbox{}\hskip2cm
\times
\left(\sum_{k=k_0}^{k'}\sum_{j=j_0}^{J_0}\frac{||D^{\sigma} P_j u||_{k}^{k}}{2^{\hat B_{k}}}
+\sum_{k=k_0}^{k'}\sum_{j=j_0}^{J_0}
\left(\frac{||D^{\sigma} P_j u||_{k}^{k}}{2^{\hat B_{k}}}\right)^{5}\right)
\end{array}\eqno{(6.10)}$$
for any $0\leq t\leq T$, where $C$ is a universal constant.}

\proof
Let $2k_0\leq k\leq k'$, $1\leq j\leq J_0$ and $0\leq t\leq T$.
We divide the proof into five steps.

\bigskip

{\it  Step 1.}\hsa
By the same arguments to derive (6.7), we have
$$\begin{array}{l}\displaystyle\mbox{}\hskip0cm
\frac{d}{dt}||D^{\sigma} P_j u||_{k}^{k}+
\frac{\nu}4k^2
\int_{R^3}|D^{\sigma} P_j u|^{k-2}|\nabla|D^{\sigma}P_j u||^2
\\[15pt]\displaystyle\mbox{}\hskip2cm
\leq
\frac{ Ck^4}{\nu}||D^{\sigma} P_j\left(u\otimes u\right)||_k^2
||D^{\sigma} P_ju||_k^{k-2},
\end{array}\eqno{(6.11)}$$
where $k(k-2)\geq k^2/2$ is used.

\bigskip

{\it Step 2. Estimate of $
||D^{\sigma} P_j\left(u\otimes u\right)||_k$ in (6.11).}\hsa
As in Step 3 of the proof of Lemma 6.2, by Lemma 2.1, 2.3 and 2.6, we have
$$\begin{array}{l}\displaystyle\mbox{}\hskip0cm
||D^{\sigma} P_j\left(u\otimes u\right)||_k
\leq C2^{(\sigma+3-\frac3k)j}||P_j\left(u\otimes u\right)||_1
\leq C2^{(\sigma+3)J_0}||u||_2^2,
\end{array}$$
where $j\leq J_0$ is used.

As $k\leq\frac{2^B}{\frac{8(\sigma+3)}{\sigma}B}$, from (6.2) and (6.3), we have
$$\frac{2^{(\sigma+3)J_0}}{2^{\hat B_k/k}}\leq
\frac{2^{\frac{8(\sigma+3)}{\sigma}B}}{2^{\frac{2^B}{k}}}\leq1\leq2^{\frac{\sqrt{k}}{1000}},
$$
for any $B>1$.

As $k>\frac{2^B}{\frac{8(\sigma+3)}{\sigma}B}$, we have
$$\frac{2^{(\sigma+3)J_0}}{2^{\hat B_k/k}}\leq2^{(\sigma+3)J_0}\leq2^{\frac{8(\sigma+3)}{\sigma}B}
\leq C2^{\frac1{1000}\sqrt{\frac{2^B}{\frac{8(\sigma+3)}{\sigma}B}}}
\leq C2^{\frac{\sqrt{k}}{1000}}$$
for any $B>1$, where $C$ is a universal constant.

Therefore
$$2^{(\sigma+3)J_0}\leq C2^{\hat B_k/k}2^{\frac{\sqrt{k}}{1000}}$$
and then
$$\begin{array}{l}\displaystyle\mbox{}\hskip0cm
||D^{\sigma} P_j\left(u\otimes u\right)||_k
\leq
C2^{\hat B_k/k}2^{\frac{\sqrt{k}}{1000}}||u||_2^2
\end{array}$$
for any $B>1$ and $k\geq2k_0$.

Plug the above inequality into (6.11) and we conclude
$$\begin{array}{l}\displaystyle\mbox{}\hskip0cm
\frac{d}{dt}||D^{\sigma} P_j u||_{k}^{k}
+
\frac{\nu}4k^2
\int_{R^3}|D^{\sigma} P_j u|^{k-2}|\nabla|D^{\sigma}P_j u||^2
\\[15pt]\displaystyle\mbox{}\hskip2cm
\leq
\frac{ Ck^4}{\nu}2^{2\hat B_k/k}2^{\frac{\sqrt{k}}{100}}||u||_2^4||D^{\sigma} P_ju||_k^{k-2}.
\end{array}\eqno{(6.12)}$$

\bigskip

{\it Step 3. Gain from $\frac{\nu}4k^2
\int_{R^3}|D^{\sigma} P_j u|^{k-2}|\nabla|D^{\sigma}P_j u||^2$.}\hsa
From interpolation inequality, we deduce,
$$||D^{\sigma} P_j u||_{k}\leq\left(||D^{\sigma} P_j u||_{\left[\frac{k+1}2\right]}\right)^{1-\lambda}
\left(||D^{\sigma} P_j u||_{3k}\right)^\lambda,$$
where $$\lambda=\frac{\frac1{\left[\frac{k+1}2\right]}-\frac1k}
{\frac1{\left[\frac{k+1}2\right]}-\frac1{3k}}.$$
It is clear that
$$\left\{\begin{array}{l}\displaystyle\mbox{}\hskip0cm
\lambda\geq\frac{\frac2{k+1}-\frac1k}
{\frac2{k+1}-\frac1{3k}}=\frac{3k-3}{5k-1}\geq\frac25\mb{and}
\\[15pt]\displaystyle\mbox{}\hskip0cm
\lambda\leq\frac{\frac2{k}-\frac1k}{\frac2{k}-\frac1{3k}}=\frac35.
\end{array}\right.\eqno{(6.13)}$$
By Lemma 4.1, there exists a universal constant $C$ such that
$$||D^{\sigma} P_j u||_{3k}^k\leq Ck^2\int_{R^3}|D^{\sigma} P_j u|^{k-2}|\nabla|D^{\sigma}P_j u||^2.$$
For simplicity we denote $\frac{ Ck^4}{\nu}2^{2\hat B_k/k}2^{\frac{\sqrt{k}}{100}}||u||_2^4$ by $\Theta$.
It follows that
$$\begin{array}{l}\displaystyle\mbox{}\hskip0cm
\mbox{Righthand side of}\hsa (6.12)=\Theta||D^{\sigma} P_ju||_k^{k-2}
\leq\Theta
\left(||D^{\sigma} P_j u||_{\left[\frac{k+1}2\right]}^{k-2}\right)^{1-\lambda}
\left(||D^{\sigma} P_j u||_{3k}^{k-2}\right)^\lambda
\\[15pt]\displaystyle\mbox{}\hskip0cm
\leq
\left(\Theta^{\frac1{1-\lambda}}||D^{\sigma} P_j u||_{\left[\frac{k+1}2\right]}^{k-2}\right)^{1-\lambda}
\left(Ck^2\int_{R^3}|D^{\sigma} P_j u|^{k-2}|\nabla|D^{\sigma}P_j u||^2\right)^{\frac{k-2}{k}\lambda}
\\[15pt]\displaystyle\mbox{}\hskip0cm
=
\left(\left(\frac{4C}{\nu}\right)^{\frac{\frac{k-2}{k}\lambda}{1-\lambda}}\Theta^{\frac1{1-\lambda}}||D^{\sigma} P_j u||_{\left[\frac{k+1}2\right]}^{k-2}\right)^{1-\lambda}
\left(\frac{\nu}4k^2\int_{R^3}|D^{\sigma} P_j u|^{k-2}|\nabla|D^{\sigma}P_j u||^2\right)^{\frac{k-2}{k}\lambda}
\\[15pt]\displaystyle\mbox{}\hskip0cm
\leq
\left(\frac{4C}{\nu}\right)^{\frac{\frac{k-2}{k}\lambda}{1-\frac{k-2}{k}\lambda}}
\Theta^{\frac1{1-\frac{k-2}{k}\lambda}}
||D^{\sigma} P_j u||_{\left[\frac{k+1}2\right]}^{\frac{(k-2)(1-\lambda)}{1-\frac{k-2}{k}\lambda}}+
\frac{\nu}4k^2\int_{R^3}|D^{\sigma} P_j u|^{k-2}|\nabla|D^{\sigma}P_j u||^2
\\[15pt]\displaystyle\mbox{}\hskip0cm
=
\left(\frac{4C}{\nu}\right)^{\frac{\frac{k-2}{k}\lambda}{1-\frac{k-2}{k}\lambda}}
\left(\frac{ Ck^4}{\nu}2^{\frac{\sqrt{k}}{100}}||u||_2^4\right)^{\frac1{1-\frac{k-2}{k}\lambda}}
2^{\hat B_k\frac{\frac2k}{1-\frac{k-2}{k}\lambda}}
||D^{\sigma} P_j u||_{\left[\frac{k+1}2\right]}^{k\frac{(1-\frac2k)(1-\lambda)}{1-\frac{k-2}{k}\lambda}}
\\[15pt]\displaystyle\mbox{}\hskip2cm
+
\frac{\nu}4k^2\int_{R^3}|D^{\sigma} P_j u|^{k-2}|\nabla|D^{\sigma}P_j u||^2.
\end{array}$$
By (6.13), it is easy to see that
$$\begin{array}{l}\displaystyle\mbox{}\hskip0cm
\left(\frac{4C}{\nu}\right)^{\frac{\frac{k-2}{k}\lambda}{1-\frac{k-2}{k}\lambda}}
\left(\frac{ Ck^4}{\nu}2^{\frac{\sqrt{k}}{100}}||u||_2^4\right)^{\frac1{1-\frac{k-2}{k}\lambda}}
\\[15pt]\displaystyle\mbox{}\hskip3cm
\leq
C\left(1+\frac1{\nu}\right)^{4}k^{10}2^{\frac{\sqrt{k}}{10}}\left(1+||u||_2^2\right)^5.
\end{array}$$

Plug this estimate into (6.12) and we obtain
$$\begin{array}{l}\displaystyle\mbox{}\hskip0cm
\frac{d}{dt}||D^{\sigma} P_j u||_{k}^{k}
\leq C\left(1+\frac1{\nu}\right)^{4}k^{10}2^{\frac{\sqrt{k}}{10}}\left(1+||u||_2^2\right)^5
\\[15pt]\displaystyle\mbox{}\hskip3cm\times
2^{\hat B_k\frac{\frac2k}{1-\frac{k-2}{k}\lambda}}
||D^{\sigma} P_j u||_{\left[\frac{k+1}2\right]}^{k\frac{(1-\frac2k)(1-\lambda)}{1-\frac{k-2}{k}\lambda}}.
\end{array}\eqno{(6.14)}$$

\bigskip

{\it Step 4. Dividing by $2^{\hat B_k}$.}\hsa
We divide by $2^{\hat B_k}$ on both sides of (6.14) and consequently,
$$\begin{array}{l}\displaystyle\mbox{}\hskip0cm
\frac{d}{dt}\frac{||D^{\sigma} P_j u||_{k}^{k}}{2^{\hat B_k}}
\leq
C\left(1+\frac1{\nu}\right)^{4}k^{10}2^{\frac{\sqrt{k}}{10}}\left(1+||u||_2^2\right)^5
\left(\frac{||D^{\sigma} P_j u||_{\left[\frac{k+1}2\right]}^k}{2^{\hat B_k}}\right)^{\frac{(1-\frac2k)(1-\lambda)}{1-\frac{k-2}{k}\lambda}},
\end{array}$$
where
$\frac{\frac2k}{1-\frac{k-2}{k}\lambda}+\frac{(1-\frac2k)(1-\lambda)}{1-\frac{k-2}{k}\lambda}=1$ is used.

In view of (6.2), (6.13) and $k\geq2k_0=200$,
$$\begin{array}{l}\displaystyle\mbox{}\hskip0cm
\left(\frac{2^{\hat B_{\left[(k+1)/2\right]}\frac{k}{\left[(k+1)/2\right]}}}{2^{\hat B_k}}\right)
^{\frac{(1-\frac2k)(1-\lambda)}{1-\frac{k-2}{k}\lambda}}
=\left(2^{k\left(\frac1{\sqrt{k}}-\frac{1}{\sqrt{\left[(k+1)/2\right]}}\right)}\right)
^{\frac{(1-\frac2k)(1-\lambda)}{1-\frac{k-2}{k}\lambda}}
\\[15pt]\displaystyle\mbox{}\hskip1cm
\leq2^{\frac35k\left(\frac1{\sqrt{k}}-\frac{1}{\sqrt{\frac23k}}\right)}=2^{\frac35\sqrt{k}\left(1-\sqrt{3/2}\right)}
\leq Ck^{-10}2^{-\frac{\sqrt{k}}{10}},
\end{array}$$
where $C$ is universal constant. It follows that
$$\begin{array}{l}\displaystyle\mbox{}\hskip0cm
\left(\frac{||D^{\sigma} P_j u||_{\left[\frac{k+1}2\right]}^k}{2^{\hat B_k}}\right)^{\frac{(1-\frac2k)(1-\lambda)}{1-\frac{k-2}{k}\lambda}}
\\[15pt]\displaystyle\mbox{}\hskip1cm
=\left(\frac{2^{\hat B_{\left[(k+1)/2\right]}\frac{k}{\left[(k+1)/2\right]}}}{2^{\hat B_k}}\right)
^{\frac{(1-\frac2k)(1-\lambda)}{1-\frac{k-2}{k}\lambda}}
\left(\frac{||D^{\sigma} P_j u||_{\left[\frac{k+1}2\right]}^{\left[\frac{k+1}2\right]}}{2^{\hat B_{\left[\frac{k+1}2\right]}}}\right)^{\frac{k}{\left[\frac{k+1}2\right]}\frac{(1-\frac2k)(1-\lambda)}{1-\frac{k-2}{k}\lambda}}
\\[15pt]\displaystyle\mbox{}\hskip1cm
\leq Ck^{-10}2^{-\frac{\sqrt{k}}{10}}\left(\frac{||D^{\sigma} P_j u||_{\left[\frac{k+1}2\right]}^{\left[\frac{k+1}2\right]}}{2^{\hat B_{\left[\frac{k+1}2\right]}}}\right)^{\frac{k}{\left[\frac{k+1}2\right]}\frac{(1-\frac2k)(1-\lambda)}{1-\frac{k-2}{k}\lambda}}
\\[15pt]\displaystyle\mbox{}\hskip1cm
\leq Ck^{-10}2^{-\frac{\sqrt{k}}{10}}\left(\frac{||D^{\sigma} P_j u||_{\left[\frac{k+1}2\right]}^{\left[\frac{k+1}2\right]}}{2^{\hat B_{\left[\frac{k+1}2\right]}}}
+
\left(\frac{||D^{\sigma} P_j u||_{\left[\frac{k+1}2\right]}^{\left[\frac{k+1}2\right]}}{2^{\hat B_{\left[\frac{k+1}2\right]}}}\right)^{5}\right)
,
\end{array}$$
where $1\leq\frac{k}{\left[\frac{k+1}2\right]}\frac{(1-\frac2k)(1-\lambda)}{1-\frac{k-2}{k}\lambda}\leq5$ is used.
Therefore
$$\begin{array}{l}\displaystyle\mbox{}\hskip0cm
\frac{d}{dt}\frac{||D^{\sigma} P_j u||_{k}^{k}}{2^{\hat B_k}}
\leq
C\left(1+\frac1{\nu}\right)^{4}\left(1+||u||_2^2\right)^5
\\[15pt]\displaystyle\mbox{}\hskip3cm
\times\left(\frac{||D^{\sigma} P_j u||_{\left[\frac{k+1}2\right]}^{\left[\frac{k+1}2\right]}}{2^{\hat B_{\left[\frac{k+1}2\right]}}}
+
\left(\frac{||D^{\sigma} P_j u||_{\left[\frac{k+1}2\right]}^{\left[\frac{k+1}2\right]}}{2^{\hat B_{\left[\frac{k+1}2\right]}}}\right)^{5}\right)
.
\end{array}\eqno{(6.15)}$$

\bigskip

{\it Step 5. Proof of (6.10).}\hsa
Take sum of $j$ from $j_0$ to $J_0$ and of $k$ from $2k_0$ to $k'$ on both sides of (6.15) and consequently,
$$\begin{array}{l}\displaystyle\mbox{}\hskip0cm
\frac{d}{dt}\sum_{k=2k_0}^{k'}\sum_{j=j_0}^{J_0}\frac{||D^{\sigma} P_j u||_{k}^{k}}{2^{\hat B_k}}
\leq
C\left(1+\frac1{\nu}\right)^{4}\left(1+||u||_2^2\right)^5
\\[15pt]\displaystyle\mbox{}\hskip0.5cm
\times\left(\sum_{k=2k_0}^{k'}\sum_{j=j_0}^{J_0}\frac{||D^{\sigma} P_j u||_{\left[\frac{k+1}2\right]}^{\left[\frac{k+1}2\right]}}{2^{\hat B_{\left[\frac{k+1}2\right]}}}
+\sum_{k=2k_0}^{k'}\sum_{j=j_0}^{J_0}
\left(\frac{||D^{\sigma} P_j u||_{\left[\frac{k+1}2\right]}^{\left[\frac{k+1}2\right]}}{2^{\hat B_{\left[\frac{k+1}2\right]}}}\right)^{5}\right).
\end{array}$$
Then (6.10) follows clearly.
$\Box$

\bigskip

\pf{Theorem 6.1} Let $k'\geq2k_0$. Adding (6.5) and (6.10) together, we obtain
$$\begin{array}{l}\displaystyle\mbox{}\hskip-0cm\frac{d}{dt}\sum_{k=k_0}^{k'}\sum_{j=j_0}^{J_0}
\frac{||D^{\sigma}P_ju||_{k}^{k}}{2^{\hat B_k}}
\leq

\frac{ C}{\nu}||u||_2^4\left(2^{-\frac14B}+
\sum_{k=k_0}^{2k_0-1}\sum_{j=j_0}^{J_0}\frac{||D^{\sigma}P_ju||_{k}^{k}}{2^{\hat B_{k}}}\right)

\\[15pt]\displaystyle\mbox{}\hskip0cm
+
C\left(1+\frac1{\nu}\right)^{4}
\left(1+||u||_2^2\right)^5
\left(\sum_{k=k_0}^{k'}\sum_{j=j_0}^{J_0}\frac{||D^{\sigma} P_j u||_{k}^{k}}{2^{\hat B_{k}}}
+\sum_{k=k_0}^{k'}\sum_{j=j_0}^{J_0}
\left(\frac{||D^{\sigma} P_j u||_{k}^{k}}{2^{\hat B_{k}}}\right)^{5}\right).
\end{array}$$
It follows that
$$\begin{array}{l}\displaystyle\mbox{}\hskip-0cm\frac{d}{dt}\sum_{k=k_0}^{k'}\sum_{j=j_0}^{J_0}
\frac{||D^{\sigma}P_ju||_{k}^{k}}{2^{\hat B_k}}
\leq
C\left(1+\frac1{\nu}\right)^{4}
\left(1+||u||_2^2\right)^5
\\[15pt]\displaystyle\mbox{}\hskip1cm\times
\left(2^{-\frac{B}4}+\sum_{k=k_0}^{k'}\sum_{j=j_0}^{J_0}\frac{||D^{\sigma} P_j u||_{k}^{k}}{2^{\hat B_{k}}}
+\sum_{k=k_0}^{k'}\sum_{j=j_0}^{J_0}
\left(\frac{||D^{\sigma} P_j u||_{k}^{k}}{2^{\hat B_{k}}}\right)^{5}\right).
\end{array}\eqno{(6.16)}$$

Set
$$F(t)=\sum_{k=k_0}^{k'}\sum_{j=j_0}^{J_0}\frac{||D^{\sigma} P_j u||_{k}^{k}}{2^{\hat B_k}}$$
and
$$g(t)=C\left(1+\frac1{\nu}\right)^{4}\left(1+||u||_2^2\right)^5.$$
Then (6.16) implies
$$\begin{array}{l}\displaystyle\mbox{}\hskip-0cm\frac{d}{dt}F(t)
\leq g(t)\left(2^{-\frac{B}4}+F(t)+F^5(t)\right),
\end{array}\eqno{(6.17)}$$
where $$\sum_{k=k_0}^{k'}\sum_{j=j_0}^{J_0}
\left(\frac{||D^{\sigma} P_j u||_{k}^{k}}{2^{\hat B_{k}}}\right)^{5}\leq\left(\sum_{k=k_0}^{k'}\sum_{j=j_0}^{J_0}
\frac{||D^{\sigma} P_j u||_{k}^{k}}{2^{\hat B_{k}}}\right)^{5}$$ is used.
From Theorem 3.1, it is easy to see that
$$\int_0^Tg(t)dt\leq C\left(1+\frac1{\nu}\right)^{4}\left(1+||u_0||_2^2\right)^5T.$$
Therefore there exists $B_1'>0$ such that if $B\geq B_1'$,
$$2^{-\frac{B}4}\leq\frac{1}{\left(3e^{\int_0^Tg(t)dt}\right)^{\frac54}}\eqno{(6.18)}$$
or
$$B_1'\geq C\left(1+\frac1{\nu}\right)^{4}\left(1+||u_0||_2^2\right)^5T$$
for some universal constant $C$.

From Lemma 4.7, there exists $\tilde B_0>0$ such that if $$B\geq\tilde B_0,$$ then
$$F(0)=\sum_{k=k_0}^{k'}\sum_{j=j_0}^{J_0}\frac{||D^{\sigma} P_j u_0||_{k}^{k}}{2^{\hat B_k}}\leq\sum_{k=k_0}^{\infty}\sum_{j=j_0}^{\infty}\frac{||D^{\sigma}P_j
u_0||_{k}^{k}}{2^{B(1-\frac{1}{\sqrt{k}})k}}\leq2^{-\frac{B}4},$$
where $\hat B_k\geq B(1-\frac{1}{\sqrt{k}})k$ is used.

Set $$\tilde B_1=\max\{B_1',\tilde B_0\}.$$
If $B\geq\tilde B_1$, from (6.17) and Lemma 4.3, we have
$$F(t)\leq2^{-\frac{B}4}3e^{\int_0^Tg(t)dt}\leq1$$
for any $0\leq t\leq T$, where (6.18) is used. That is,
$$\sum_{k=k_0}^{k'}\sum_{j=j_0}^{J_0}\frac{||D^{\sigma} P_j u||_{k}^{k}}{2^{\hat B_k}}\leq1$$
for any $k'\geq2k_0$. This implies (6.4) clearly.
$\Box$

\bigskip

\section{Regularity improving}

In this section, we will prove the following Theorem 7.1 from which we can show the convergence of the high frequency part (large $j$) of series (5.1). This kind of regularity improving is not new essentially, but here we need a special form.

\bigskip

\theorem{7.1}{
Let $k_0=100$, $j_0=1$, $\sigma=2$, $0\leq T'\leq T$, $\cal B$ and $B>1$ be real numbers, $B_k$ be defined by (5.2) for any $k\geq k_0$ and $u_0$ be a function satisfying (1.2).
Suppose $u$ and $p$ satisfy Condition (S). There exists $\tilde B_2>1$ depending only on $u_0$ and $\sigma$ such that
if $B\geq\tilde B_2$ and
$$\displaystyle
\sum_{j=j_0}^{\infty}||D^{\sigma}P_{j}u(t)||_{k_0}^{k_0}\leq{\cal B}2^{B_{k_0}}\eqno{(7.1)}$$
for any $0\leq t\leq T'$, then we have
$$\begin{array}{l}\displaystyle
\sum_{j=j_0}^{\infty}||D^{\sigma+1}P_ju(t)||_{k_0}^{k_0}
\leq C\left(1+\frac1{\nu}\right)^{k_0+1}\left(1+T\right)
\left(1+||u_0||_2\right)^{2k_0+2}
\\[15pt]\displaystyle\mbox{}\hskip1cm\times
{\cal B}^{1+\frac{5(k_0+2)/2}{(\sigma+3/2)k_0-3}}
2^{\frac{5(k_0+2)B/2}{\sigma+3/2-3/k_0}}
2^{B_{k_0}}
\end{array}\eqno{(7.2)}$$
for any $0\leq t\leq T'$, where $C>0$ is a universal constant.

}

\bigskip

\proof We divide the proof into four steps.

\bigskip

{\it Step 1. Estimate of $||u||_{\infty}$.}\hsa(7.1) implies
$$\displaystyle
||D^{\sigma}P_{j}u(t)||_{k_0}\leq\left({\cal B}2^{B_{k_0}}\right)^{\frac1{k_0}}
\leq4{\cal B}^{\frac1{k_0}}2^B$$
for any $j\geq j_0$ and $0\leq t\leq T'$, where (5.2) is used.
From Lemma 2.1, 2.3 and 2.6, we have
$$\displaystyle
||DP_ju||_{\infty}\leq C2^{(\frac3{k_0}+1-\sigma)j}||D^{\sigma}P_ju||_{k_0}$$
and
$$
||DP_ju||_{\infty}\leq C2^{\frac52j}||P_ju||_{2}.
$$
It follows that
$$\begin{array}{l}\displaystyle
||DP_ju||_{\infty}\leq C||P_ju||_{2}^{\frac{\sigma-1-3/k_0}{\sigma+3/2-3/k_0}}
||D^{\sigma}P_ju||_{k_0}^{\frac{5/2}{\sigma+3/2-3/k_0}}
\\[15pt]\displaystyle\mbox{}\hskip1cm
\leq C||u_0||_{2}^{\frac{\sigma-1-3/k_0}{\sigma+3/2-3/k_0}}
{\cal B}^{\frac{5/2}{(\sigma+3/2)k_0-3}}2^{\frac{5B/2}{\sigma+3/2-3/k_0}}.
\end{array}$$

From (2.1), Lemma 2.1 and Theorem 3.1, we have
$$\begin{array}{l}\displaystyle
||u||_{\infty}\leq||P_{\leq0}u||_{\infty}+\sum_{j=j_0}^{\infty}||P_ju||_{\infty}
\leq C||u||_2+C\sum_{j=j_0}^{\infty}2^{-j}||DP_ju||_{\infty}
\\[15pt]\displaystyle\mbox{}\hskip0.5cm
\leq C||u_0||_2+C||u_0||_{2}^{\frac{\sigma-1-3/k_0}{\sigma+3/2-3/k_0}}
{\cal B}^{\frac{5/2}{(\sigma+3/2)k_0-3}}2^{\frac{5B/2}{\sigma+3/2-3/k_0}}
\\[15pt]\displaystyle\mbox{}\hskip0.5cm
\leq C\left(1+||u_0||_2\right){\cal B}^{\frac{5/2}{(\sigma+3/2)k_0-3}}2^{\frac{5B/2}{\sigma+3/2-3/k_0}},
\end{array}\eqno{(7.3)}$$
where $C$ is a universal constant.
For simplicity, we denote
$$A=\left(1+||u_0||_2\right){\cal B}^{\frac{5/2}{(\sigma+3/2)k_0-3}}2^{\frac{5B/2}{\sigma+3/2-3/k_0}}.\eqno{(7.4)}$$

\bigskip

{\it Step 2.}\hsa
By the same arguments to derive (6.11), we have
$$\begin{array}{l}\displaystyle\mbox{}\hskip0cm
\frac{d}{dt}||D^{\sigma+1} P_j u||_{k_0}^{k_0}+
\frac{\nu}4k_0^2\int_{R^3}|D^{\sigma+1} P_ju|^{k_0-2}|\nabla D^{\sigma+1}P_ju|^2
\\[15pt]\displaystyle\mbox{}\hskip2cm
\leq
\frac{ C}{\nu}||D^{\sigma+1} P_j\left(u\otimes u\right)||_{k_0}^2
||D^{\sigma+1} P_ju||_{k_0}^{k_0-2}.
\end{array}\eqno{(7.5)}$$

By Corollary 2.8, where we set $\hat q=2$, $q=q_1=k_0$ and $q_0=\infty$, we deduce
$$\begin{array}{l}\displaystyle\mbox{}\hskip0cm
||D^{\sigma+1} P_j\left(u\otimes u\right)||_{k_0}^2
\leq C\Bigg\{\alpha_j||u||_2^4+
\sum_{m=\max\{j_0,j-2\}}^{\infty}2^{(\sigma+1)(j-m)}||D^{\sigma+1}P_mu||_{k_0}^2||u||_{\infty}^2
\Bigg\},
\end{array}$$
where $\alpha_j$ is given by (2.3). In view of (7.3),
$$\begin{array}{l}\displaystyle\mbox{}\hskip0cm
||D^{\sigma+1} P_j\left(u\otimes u\right)||_{k_0}^2
\leq C\alpha_j||u||_2^4+
CA^2\sum_{m=\max\{j_0,j-2\}}^{\infty}2^{(\sigma+1)(j-m)}||D^{\sigma+1}P_mu||_{k_0}^2
.
\end{array}$$

Plug this into (7.5) and we obtain
$$\begin{array}{l}\displaystyle\mbox{}\hskip0cm
\frac{d}{dt}||D^{\sigma+1} P_j u||_{k_0}^{k_0}+
\frac{\nu}4k_0^2\int_{R^3}|D^{\sigma+1} P_ju|^{k_0-2}|\nabla D^{\sigma+1}P_ju|^2
\\[15pt]\displaystyle\mbox{}\hskip0.5cm
\leq\frac{C}{\nu}\Bigg\{\alpha_j
||u||_2^4||D^{\sigma+1} P_ju||_{k_0}^{k_0-2}
\\[15pt]\displaystyle\mbox{}\hskip1cm
+A^2\sum_{m=\max\{j_0,j-2\}}^{\infty}2^{(\sigma+1)(j-m)}||D^{\sigma+1}P_mu||_{k_0}^2
||D^{\sigma+1} P_ju||_{k_0}^{k_0-2}\Bigg\}.
\end{array}\eqno{(7.6)}$$

\bigskip

{\it Step 3. Taking sum of $j$.}\hsa
It is easy to see that
$$\begin{array}{l}\displaystyle\mbox{}\hskip0cm
\alpha_j||u||_2^4||D^{\sigma+1} P_ju||_{k_0}^{k_0-2}
\leq\alpha_j\left(||u||_2^4\right)^{\frac{k_0}2}
+||D^{\sigma+1} P_ju||_{k_0}^{k_0}
\end{array}$$
and
$$\begin{array}{l}\displaystyle\mbox{}\hskip0cm
\sum_{m=\max\{j_0,j-2\}}^{\infty}2^{(\sigma+1)(j-m)}||D^{\sigma+1}P_mu||_{k_0}^2||D^{\sigma+1} P_ju||_{k_0}^{k-2}
\\[15pt]\displaystyle\mbox{}\hskip1.6cm
\leq\sum_{m=\max\{j_0,j-2\}}^{\infty}2^{(\sigma+1)(j-m)}
\left(||D^{\sigma+1}P_mu||_{k_0}^{k_0}+||D^{\sigma+1} P_ju||_{k_0}^{k_0}\right)
\\[15pt]\displaystyle\mbox{}\hskip1.6cm
\leq C||D^{\sigma+1} P_ju||_{k_0}^{k_0}+
\sum_{m=\max\{j_0,j-2\}}^{\infty}2^{(\sigma+1)(j-m)}
||D^{\sigma+1}P_mu||_{k_0}^{k_0}.
\end{array}$$
Plug these two inequalities into (7.6) and we have
$$\begin{array}{l}\displaystyle\mbox{}\hskip0cm
\frac{d}{dt}||D^{\sigma+1} P_j u||_{k_0}^{k_0}+
\frac{\nu}4k_0^2\int_{R^3}|D^{\sigma+1} P_ju|^{k_0-2}|\nabla D^{\sigma+1}P_ju|^2
\\[15pt]\displaystyle\mbox{}\hskip0.5cm
\leq\frac{C}{\nu}\Bigg\{\alpha_j
||u||_2^{2k_0}+A^2||D^{\sigma+1} P_ju||_{k_0}^{k_0}
\\[15pt]\displaystyle\mbox{}\hskip2cm
+A^2\sum_{m=\max\{j_0,j-2\}}^{\infty}2^{(\sigma+1)(j-m)}||D^{\sigma+1}P_mu||_{k_0}^{k_0}
\Bigg\},
\end{array}$$
where $A>1$ is used.

Take sum of $j$ from $j_0=1$ to $\infty$ on both sides of the above inequality.
Since
$$\begin{array}{l}\displaystyle\mbox{}\hskip0cm
\sum_{j=1}^{\infty}\sum_{m=\max\{j_0,j-2\}}^{\infty}2^{(\sigma+1)(j-m)}||D^{\sigma+1}P_mu||_{k_0}^{k_0}
=
\sum_{m=1}^{\infty}2^{(\sigma+1)(1-m)}||D^{\sigma+1}P_mu||_{k_0}^{k_0}
\\[15pt]\displaystyle\mbox{}\hskip1cm
+
\sum_{m=1}^{\infty}2^{(\sigma+1)(2-m)}||D^{\sigma+1}P_mu||_{k_0}^{k_0}+
\sum_{j=3}^{\infty}\sum_{m=j-2}^{\infty}2^{(\sigma+1)(j-m)}||D^{\sigma+1}P_mu||_{k_0}^{k_0}
\\[15pt]\displaystyle\mbox{}\hskip0cm
\leq C\sum_{m=1}^{\infty}||D^{\sigma+1}P_mu||_{k_0}^{k_0}+
\sum_{m=1}^{\infty}\sum_{j=3}^{m+2}2^{(\sigma+1)(j-m)}||D^{\sigma+1}P_mu||_{k_0}^{k_0}
\\[15pt]\displaystyle\mbox{}\hskip0cm
\leq
C\sum_{j=1}^{\infty}||D^{\sigma+1}P_ju||_{k_0}^{k_0},
\end{array}$$
using Lemma 4.2 and 4.6, we have
$$\begin{array}{l}\displaystyle\mbox{}\hskip0cm
\frac{d}{dt}\sum_{j=1}^{\infty}||D^{\sigma+1} P_j u||_{k_0}^{k_0}+
\frac{\nu}4k_0^2\sum_{j=1}^{\infty}\int_{R^3}|D^{\sigma+1} P_ju|^{k_0-2}|\nabla D^{\sigma+1}P_ju|^2
\\[15pt]\displaystyle\mbox{}\hskip1.5cm
\leq
\frac{C}{\nu}\left(||u||_2^{2k_0}
+
A^2\sum_{j=1}^{\infty}||D^{\sigma+1} P_ju||_{k_0}^{k_0}\right)
.
\end{array}\eqno{(7.7)}$$

\bigskip

{\it Step 4. Proof of (7.2).}\hsa
From Lemma 4.4, we have
$$\begin{array}{l}\displaystyle\mbox{}\hskip0cm
||D^{\sigma+1} P_ju||_{k_0}^{k_0}
\leq
C{k_0}^2\left(\int_{R^3}|D^{\sigma+1}P_ju|^{k_0-2}|\nabla D^{\sigma+1}P_ju|^2\right)^{\frac{k_0}{k_0+2}}
\left(||D^{\sigma} P_ju||_{k_0}^{k_0}\right)^{\frac{2}{k_0+2}}
\\[15pt]\displaystyle\mbox{}\hskip0.4cm
\leq\frac{\nu^2k_0^2}{4CA^2}
\int_{R^3}|D^{\sigma+1}P_ju|^{k_0-2}|\nabla D^{\sigma+1}P_ju|^2+\left(Ck_0^2\right)^{\frac{k_0+2}2}
\left(\frac{4CA^2}{\nu^2k_0^2}\right)^{\frac{k_0}2}||D^{\sigma} P_ju||_{k_0}^{k_0}
\end{array}$$
for any $0\leq t\leq T'$.
Plug this inequality into (7.7) and we obtain
$$\begin{array}{l}\displaystyle\mbox{}\hskip0cm
\frac{d}{dt}\sum_{j=1}^{\infty}||D^{\sigma+1} P_j u||_{k_0}^{k_0}
\leq
\frac{C}{\nu}||u||_2^{2k_0}
+
CA^{2+k_0}\left(\frac1{\nu}\right)^{k_0+1}\sum_{j=1}^{\infty}
||D^{\sigma} P_ju||_{k_0}^{k_0}
\\[15pt]\displaystyle\mbox{}\hskip1cm
\leq
\frac{C}{\nu}||u||_2^{2k_0}
+
C\left(\frac1{\nu}\right)^{k_0+1}A^{2+k_0}{\cal B}2^{B_{k_0}}
,
\end{array}$$
where (7.1) is used.
It follows that
$$\begin{array}{l}\displaystyle\mbox{}\hskip0cm
\sum_{j=1}^{\infty}||D^{\sigma+1} P_j u(t)||_{k_0}^{k_0}
\leq \sum_{j=1}^{\infty}||D^{\sigma+1} P_j u_0||_{k_0}^{k_0}+
\frac{C}{\nu}||u||_2^{2k_0}T'
\\[15pt]\displaystyle\mbox{}\hskip5cm
+
C\left(\frac1{\nu}\right)^{k_0+1}A^{2+k_0}{\cal B}2^{B_{k_0}}T'
\end{array}\eqno{(7.8)}$$
for any $0\leq t\leq T'$.

From Lemma 4.7, there exists $\tilde B_2$ such that if $B\geq\tilde B_2$, then
$$\sum_{j=1}^{\infty}\frac{||D^{\sigma+1} P_j u_0||_{k_0}^{k_0}}
{2^{B(1-\frac1{\sqrt{k}})k}}\leq1$$
or
$$\sum_{j=1}^{\infty}||D^{\sigma+1} P_j u_0||_{k_0}^{k_0}\leq
2^{B(1-\frac1{\sqrt{k}})k}\leq2^{B_{k_0}}.$$
Plug this into (7.8) and we arrive at
$$\begin{array}{l}\displaystyle\mbox{}\hskip0cm
\sum_{j=1}^{\infty}||D^{\sigma+1} P_j u(t)||_{k_0}^{k_0}
\leq C\left(1+\frac1{\nu}\right)^{k_0+1}\left(1+T\right)
\left(1+||u||_2^2\right)^{k_0}
A^{2+k_0}{\cal B}2^{B_{k_0}}.
\end{array}$$
In view of (7.4), we have (7.2) clearly.
$\Box$

\bigskip

\corollary{7.2}{Suppose all the assumptions of Theorem 7.1 hold.
Then if $B\geq\tilde B_2$ given by Theorem 7.1, we have
$$\begin{array}{l}\displaystyle
||D^{\sigma}P_ju(t)||_{k}
\leq\hat{\cal B}2^{\left(\frac3{k_0}-\frac3k-1\right)j}2^{\left(1+\frac5{\sigma}\right)B}
\end{array}\eqno{(7.9)}$$
for any $j\geq j_0$, $k\geq k_0$ and $0\leq t\leq T'$, where
$$\hat{\cal B}:=C\left(1+\frac1{\nu}\right)^{1+\frac1{k_0}}\left(1+T\right)^{\frac1{k_0}}
\left(1+||u_0||_2\right)^{2+\frac2{k_0}}{\cal B}^{\frac2{k_0}}\eqno{(7.10)}$$
with $C>0$ being a universal constant.
}

\proof
(7.2) implies that
$$\begin{array}{l}\displaystyle
||D^{\sigma+1}P_ju(t)||_{k_0}^{k_0}
\leq C\left(1+\frac1{\nu}\right)^{k_0+1}\left(1+T\right)
\left(1+||u_0||_2\right)^{2k_0+2}
\\[15pt]\displaystyle\mbox{}\hskip1cm\times
{\cal B}^{1+\frac{5(k_0+2)/2}{(\sigma+3/2)k_0-3}}
2^{\frac{5(k_0+2)B/2}{\sigma+3/2-3/k_0}}
2^{B_{k_0}}
\end{array}$$
or
$$\begin{array}{l}\displaystyle
||D^{\sigma+1}P_ju(t)||_{k_0}
\leq C\left(1+\frac1{\nu}\right)^{1+\frac1{k_0}}\left(1+T\right)^{\frac1{k_0}}
\left(1+||u_0||_2\right)^{2+\frac2{k_0}}
{\cal B}^{\frac2{k_0}}
2^{\frac5{\sigma}B}
2^{B_{k_0}/k_0}
\end{array}$$
for any $j\geq j_0$ and $0\leq t\leq T'$. Recall (5.2), that is $$B_{k_0}=\left(B+1+\frac1{\sqrt{k_0}}\right)k_0,$$
we see
$$\begin{array}{l}\displaystyle
||D^{\sigma+1}P_ju(t)||_{k_0}
\leq C\left(1+\frac1{\nu}\right)^{1+\frac1{k_0}}\left(1+T\right)^{\frac1{k_0}}
\left(1+||u_0||_2\right)^{2+\frac2{k_0}}
{\cal B}^{\frac2{k_0}}
2^{\left(1+\frac5{\sigma}\right)B}.
\end{array}\eqno{(7.11)}$$
From lemma 2.3 and 2.6, we have
$$\begin{array}{l}\displaystyle
||D^{\sigma}P_ju(t)||_{k}\leq C2^{\left(\frac3{k_0}-\frac3k\right)j}
||D^{\sigma}P_ju(t)||_{k_0}\leq C2^{\left(\frac3{k_0}-\frac3k-1\right)j}||D^{\sigma+1}P_ju(t)||_{k_0}
\end{array}$$
for any $k\geq k_0$ and $0\leq t\leq T'$. Combining it with (7.11), we have (7.9).
$\Box$

\bigskip

\section{New a priori estimates}

In this section, we will demonstrate our new a priori estimates of Navier-Stokes equations, Theorem 8.1 and its corollary.
Actually, to obtain Theorem 8.1, we only need to delete the condition (5.3) in Theorem 5.1 ({\it uniform bound estimate})
by choosing a suitable large $B$ in (5.2).

\bigskip

\theorem{8.1}{Let $k_0=100$, $j_0=1$, $\sigma=2$, $T>0$ and $B>0$ be real numbers,
$B_k$ be given by (5.1) for any $k\geq k_0$ and $u_0$ be a function satisfying (1.2).
Suppose $u$ and $p$ satisfy Condition (S).
Then there exists $\tilde B>0$ depending only on $T, \nu$ and $u_0$ such that if $B\geq\tilde B$, then
$$\begin{array}{l}\displaystyle
\sum_{k=k_0}^{\infty}\sum_{j=j_0}^{\infty}
\frac{||D^{\sigma}P_j
u(t)||_{k}^{k}}{2^{B_k}}
\leq2{\cal C}-1
\end{array}\eqno{(8.1)}$$
for any $0\leq t\leq T$, where
${\cal C}$ is given by (5.5).}

\bigskip

We first prove the following simple lemma.

\bigskip

\lemma{8.2}{Suppose all the assumptions of Theorem 8.1 hold and $0<T'\leq T$.
If (8.1) holds for any $0<t<T'$, then it holds
for $t=T'$.
}

\proof
For any $0\leq t<T'$, $k'\geq k_0$ and $j'\geq j_0$, from (8.1), it follows that
$$\begin{array}{l}\displaystyle
\sum_{k=k_0}^{k'}\sum_{j=j_0}^{j'}
\frac{||\dd{\sigma}P_j
u(t)||_{k}^{k}}{2^{B_{k}}}
\leq
2{\cal C}-1.
\end{array}$$
Let $t\rightarrow T'$ and using
the continuity of $||\dd{\sigma}P_ju(t)||_{k}$ with respect to $t$, we obtain
$$\begin{array}{l}\displaystyle
\sum_{k=k_0}^{k'}\sum_{j=j_0}^{j'}
\frac{||\dd{\sigma}P_j
u(T')||_{k}^{k}}{2^{B_{k}}}
\leq
2{\cal C}-1.
\end{array}$$
Since $k'$ and $j'$ are arbitrary, we have (8.1) holds as $t=T'$.
$\Box$

\bigskip

\pf{Theorem 8.1}We divide the proof into six steps.

\bigskip

{\it Step 1. Choosing of $\tilde B$.}\hsa In Theorem 7.1 (and Corollary 7.2), set
$${\cal B}=2{\cal C}.\eqno{(8.2)}$$
Let $\tilde B_3$ be large enough such that
$$\hat{\cal B}\leq2^{\frac{\tilde B_3}{\sigma}},\eqno{(8.3)}$$
where $\hat{\cal B}$ is given by (7.10).
Now we set $$\tilde B=\max\left\{\tilde B_0,\tilde B_1,\tilde B_2,\tilde B_3,\sigma+1\right\},$$
where $\tilde B_0$ is given by Lemma 4.7, $\tilde B_1$ is given by Theorem 6.1 and $\tilde B_2$ is given by Theorem 7.1.

\bigskip

{\it Step 2.}\hsa
Let $B\geq\tilde B$ and
$$T':=\sup\left\{\hat T:\mbox{(8.1) holds for any}\hsa0\leq t\leq\hat T\right\}.\eqno{(8.4)}$$
In view of Lemma 8.2 and Lemma 4.7, we have (8.1) holds as $t=T'$.

If $T'=T$, then Theorem 8.1 is true.

Next, we suppose $T'<T$ and we will derive a contradiction by it.

\bigskip

{\it Step 3.}\hsa
From Condition (S) and Lemma 4.6, we have
$$\begin{array}{l}\displaystyle
\sum_{j=j_0}^{\infty}||D^{\sigma}P_ju(t)||_{k_0}^{k_0}
\end{array}$$
is continuous. In view of (8.1) and (8.4), 
$$\begin{array}{l}\displaystyle
\sum_{j=j_0}^{\infty}||D^{\sigma}P_ju(t)||_{k_0}^{k_0}\leq\left(2{\cal C}-1\right)2^{B_{k_0}}
\end{array}$$
for any $0\leq t\leq T'$. Then there exists $\delta>0$ such that $T'+\delta\leq T$ and
$$\begin{array}{l}\displaystyle
\sum_{j=j_0}^{\infty}||D^{\sigma}P_ju(t)||_{k_0}^{k_0}\leq2{\cal C}2^{B_{k_0}}={\cal B}2^{B_{k_0}}
\end{array}$$
for any $0\leq t\leq T'+\delta$, where (8.2) is used.
Therefore by Corollary 7.2 (recall $B\geq\tilde B_2$), we have (7.9) holds.
In view of (8.3) and $B\geq\tilde B_3$, (7.9) implies
$$\begin{array}{l}\displaystyle
||D^{\sigma}P_ju(t)||_{k}
\leq2^{\left(\frac3{k_0}-\frac3k-1\right)j}2^{\left(1+\frac6{\sigma}\right)B}.
\end{array}\eqno{(8.5)}$$
for any $0\leq t\leq T'+\delta$.

\bigskip

{\it Step 4. Convergence of the high frequency part.}\hsa Let $$J_0=\left[\frac{8B}{\sigma}\right].$$
Then as $j>J_0$, from $B\geq\tilde B\geq\sigma+1$, we have
$$\begin{array}{l}\displaystyle
\left(\frac3{k_0}-1\right)j+\frac{6B}{\sigma}\leq
\left(\frac3{k_0}-1\right)\left(\frac{8B}{\sigma}-1\right)+\frac{6B}{\sigma}
\\[15pt]\displaystyle\mbox{}\hskip3cm
\leq-\frac{B}{\sigma}-\left(\frac3{k_0}-1\right)\leq0.
\end{array}$$
Combining it with (8.5), we have
$$\begin{array}{l}\displaystyle
\frac{||D^{\sigma}P_ju(t)||_{k}^k}{2^{B_k}}\leq
\left(\frac{2^{\left(\frac3{k_0}-\frac3k-1\right)j}2^{\left(1+\frac6{\sigma}\right)B}}{2^{B+1+\frac{1}{\sqrt{k}}}}\right)^k
\leq
\left(\frac{2^{\left(\frac3{k_0}-\frac3k-1\right)j}2^{\frac{6B}{\sigma}}}{2}\right)^k
\\[15pt]\displaystyle\mbox{}\hskip2cm
=2^{-3j-k}
\left(2^{\left(\frac3{k_0}-1\right)j+\frac{6B}{\sigma}}\right)^k
\leq2^{-3j-k}
\end{array}$$
for any $0\leq t\leq T'+\delta$.
Therefore
$$\begin{array}{l}\displaystyle
\sum_{k=k_0}^{\infty}\sum_{j=J_0+1}^{\infty}
\frac{||D^{\sigma}P_j
u(t)||_{k}^{k}}{2^{B_k}}
\leq\sum_{k=k_0}^{\infty}\sum_{j=J_0+1}^{\infty}2^{-3j-k}\leq1
\end{array}\eqno{(8.6)}$$
for any $0\leq t\leq T'+\delta$.

\bigskip

{\it Step 5. Convergence of the low frequency part.}\hsa
Since $B\geq\tilde B\geq\tilde B_1$, by Theorem 6.1,
we have
$$\begin{array}{l}\displaystyle
\sum_{k=k_0}^{\infty}\sum_{j=j_0}^{J_0}
\frac{||D^{\sigma}P_j
u(t)||_{k}^{k}}{2^{\hat B_k}}
\leq1
\end{array}\eqno{(8.7)}$$
for any $0\leq t\leq T$.
From (5.2) and (6.2), we have
$$\lim_{k\rightarrow\infty}\frac{\hat B_k}{B_k}=\frac{B}{B+1}.$$
Therefore there exists $\hat k\geq k_0$ such that
$$\hat B_k\leq B_k.$$
In view of (8.7),
$$\begin{array}{l}\displaystyle
\sum_{k=\hat k}^{\infty}\sum_{j=j_0}^{J_0}
\frac{||D^{\sigma}P_j
u(t)||_{k}^{k}}{2^{B_k}}\leq
\sum_{k=\hat k}^{\infty}\sum_{j=j_0}^{J_0}
\frac{||D^{\sigma}P_j
u(t)||_{k}^{k}}{2^{\hat B_k}}
\leq1
\end{array}$$
for any $0\leq t\leq T$.
From condition (S), there exists $\tilde{\cal B}$ such that
$$\begin{array}{l}\displaystyle
\sum_{k=k_0}^{\hat k-1}\sum_{j=j_0}^{J_0}
\frac{||D^{\sigma}P_j
u(t)||_{k}^{k}}{2^{B_k}}\leq\tilde{\cal B}
\end{array}$$
for any $0\leq t\leq T$.
Therefore
$$\begin{array}{l}\displaystyle
\sum_{k=k_0}^{\infty}\sum_{j=j_0}^{J_0}
\frac{||D^{\sigma}P_j
u(t)||_{k}^{k}}{2^{B_k}}
\leq1+\tilde{\cal B}
\end{array}\eqno{(8.8)}$$
for any $0\leq t\leq T$.

\bigskip

{\it Step 6. Contradiction.}\hsa
From (8.6) and (8.8), we have (5.3) holds for $0\leq t\leq T'+\delta$. Then by Theorem 5.1, we have
(5.4) holds for $0\leq t\leq T'+\delta$. Since $B\geq\tilde B_0$, from Lemma 4.7, we have
$$\sum_{k=k_0}^{\infty}\sum_{j=j_0}^{\infty}\frac{||D^{\sigma}P_j
u_0||_{k}^{k}}{2^{B_k}}\leq\sum_{k=k_0}^{\infty}\sum_{j=j_0}^{\infty}\frac{||D^{\sigma}P_j
u_0||_{k}^{k}}{2^{B(1-\frac1{\sqrt{k}})k}}\leq1.$$
Then (5.4) implies (8.1) clearly. That is,
(8.1) holds for $0\leq t\leq T'+\delta$.
This contradicts with (8.4), the choice of $T'$.
$\Box$

\bigskip

\corollary{8.2}{Let $T>0$ and $u_0$ be a function satisfying (1.2).
There exists $\check B>0$ depending only on $T, \nu$ and $u_0$ such that
if $u$ and $p$ satisfy Condition (S), then
$$||u(t)||_{\infty}\leq\check B$$
for any $0\leq t\leq T$.
}

\proof Let $\tilde B$ be given by Theorem 8.1. Then we have
$$\begin{array}{l}\displaystyle
\sum_{k=k_0}^{\infty}\sum_{j=j_0}^{\infty}
\frac{||D^{\sigma}P_j
u(t)||_{k}^{k}}{2^{B_k}}
\leq2{\cal C}-1
\end{array}$$
for any $0\leq t\leq T$, where $B_k=\left(\tilde B+1+\frac1{\sqrt{k}}\right)k$ (recall (5.2)). This implies
$$\begin{array}{l}\displaystyle
||D^{\sigma}P_ju(t)||_{k}\leq\left(\left(2{\cal C}-1\right){2^{B_k}}\right)^{\frac1k}
\end{array}$$
for any $k\geq k_0$, $j\geq j_0$ and $0\leq t\leq T$. Let $k\rightarrow\infty$ and then
$$\begin{array}{l}\displaystyle
||D^{\sigma}P_ju(t)||_{\infty}\leq2^{\tilde B+1}.
\end{array}$$
It follows that
$$\begin{array}{l}\displaystyle
||u||_{\infty}=||P_{\leq0}u||_{\infty}+\sum_{j=j_0}^{\infty}||P_ju||_{\infty}
\\[15pt]\displaystyle\mbox{}\hskip1cm
\leq C||u_0||_2+C\sum_{j=j_0}^{\infty}2^{-\sigma j}||D^{\sigma}P_ju(t)||_{\infty}
\\[15pt]\displaystyle\mbox{}\hskip1cm
\leq C||u_0||_2+C2^{\tilde B+1}:=\check B.
\end{array}$$
The proof is complete.
$\Box$

\bigskip

\section{Proof of Theorem 1.1}

Theorem 1.1 is an easy consequence of our new a priori estimates.

\bigskip

\pf{Theorem 1.1}From Theorem 3.4 and 3.5, there exist $T^*>0$, smooth functions $p(x,t)$ and $u(x,t)$ on
$R^3\times[0,T^*]$ with $u(x,0)=u_0(x)$ and $f\equiv0$ such that (1.1) holds. By Theorem 3.1, we have (1.3) holds.
Let
$$\begin{array}{l}\displaystyle
T=\sup\Big\{T':\mbox{There exist smooth functions}\hsa p(x,t)\mb{and}u(x,t)\mb{on}
R^3\times[0,T']
\\[10pt]\displaystyle\mbox{}\hskip1.4cm
\mb{with}u(x,0)=u_0(x)\mb{and}f\equiv0\mb{such that (1.1) and (1.3) hold.}\Big\}.
\end{array}$$
Then $T\geq T^*$.

If $T=+\infty$, then Theorem 1.1 is true.

If $T<\infty$, then from Theorem 3.6, we have
$$\limsup_{t\rightarrow T^-}||u(t)||_{\infty}=\infty.\eqno{(9.1)}$$
Let $\check B$ be given by Corollary 8.2 with $u_0$, $\nu$ and this $T$. Then from Theorem 3.5 and Corollary 8.2, we have
$$||u(t)||_{\infty}<\check B$$
for any $0\leq t<T$.
This contradict with (9.1).

The proof of Theorem 1.1 is complete.
$\Box$

\bigskip

\remark{9.1} In Corollary 8.2, the bound of $||u||_{\infty}$, $\check B$ depends on $||D^{\sigma+3}u_0||_2$
and $||D^{\sigma+3}u_0||_{\infty}$ (recall the choosing of $\tilde B_2$ in Theorem 7.1). If $u_0\in H^1(R^3)$, for any $T>0$, we can bound
$$\sup_{0\leq t\leq T}||\nabla u(t)||_2$$
by the following way. From Theorem 3.4, there exists $T^*>0$ such that
$$\sup_{0\leq t\leq T^*}||\nabla u(t)||_2\leq\check B_1\eqno({9.2)}$$
which is a constant depending on $u_0$ and $\nu$.
By Remark 3.7, $u$ is smooth on $[T^*/2,T^*]$. Then we can use Corollary 8.2 for $t\in[T^*/2,T]$, that is,
there exists a constant $\check B_2$ depending only on $u_0$, $\nu$ and $T$ such that
$$\sup_{T^*/2\leq t\leq T}||u(t)||_{\infty}\leq\check B_2.$$
Then from the classical regularity results of parabolic equations, we have
$$\sup_{T^*/2\leq t\leq T}||\nabla u(t)||_{2}\leq\check B_3\eqno({9.3)}$$
which is a constant depending on $u_0$, $\nu$ and $T$.
From (9.2) and (9.3), we see $$\sup_{0\leq t\leq T}||\nabla u(t)||_2\leq\max\{\check B_1,\check B_2\}.\eqno{(9.4)}$$

Using the a priori estimate (9.4), if $u_0\in H^1(R^3)$,
we can conclude that (1.1) has the strong solution on $[0,T]$ for any $T>0$, which is smooth in $(0,\infty)$.
$\Box$

\bigskip

{\bf Acknowledgement.} The author would like to thank Prof.Lihe Wang
for his instruction and encouragement. The author would also like to
thank Prof.Gerhard Str\"{o}hmer and Prof.Lizhou Wang for some useful
discussions.

\bigskip

 {\Large\bf Reference}

\def\toto#1#2{\centerline{\hbox
to0.7cm{[#1]}\parbox[t]{11.5cm}{\baselineskip=12pt#2}}\vskip0.3cm}

\bigskip

{\parindent=0pt

\toto{1}{L.Caffarelli, R.V.Kohn and L.Nirenberg, Partial regularity
of suitable weak solutions of the Navier-Stokes equations, Comm.
Pure Appl. Math., Vol.35(1982), 771-831.}

\toto{2}{P.Constantin and C.Fefferman, Direction of vorticity and
the problem of global regularity for the Navier-Stokes equations,
Indiana Univ. Math. J., 42(1993), 775-789.}

\toto{3}{C.Fefferman, Existence and smoothness of the Navier-Stokes
 equation, http://www.claymath.org/millennium/Navier\_Stokes\_Equations/navierstokes.pdf.}

\toto{4}{J.Leray, Sur le mouvement d'un liquide visques emplissent
l'espace, Acta Math.J., 63(1934), 193-248.}

\toto{5}{F.H.Lin, A new proof of the Caffarelli-Kohn-Nirenberg
theorem, Comm. Pure and  Appl. math., 51(1998), 241-257.}

\toto{6}{V.Scheffer, Turbulence and Hausdorff dimension, in
Turbulence and the Navier-Stokes equations, Lecture Notes in Math.
565, Springer Verlag, Berlin, 1976, 94-112.}

\toto{7}{J.Serrin, On the interior regularity of weak solutions of
the Navier-Stokes equations, Arch. Ration. Mech. Anal., 9(1962),
187-195.}

\toto{8}{T.Tao, Lecture notes on Harmonic analysis in the phase
plane, http://www.math.ucla.edu/~tao/254a.1.01w/.}

\toto{9}{R.Teman, Navier-Stokes
 equations: Theory and numerical analysis, North-Holland-Amsterdan, New York, Oxford, 1984.}

}

\end{document}